\documentclass[11pt]{amsart}

\usepackage{amsmath,amsthm, amssymb}
\usepackage[T1]{fontenc}
\usepackage{latexsym}
\usepackage{subfigure, epsfig, epsf}
\usepackage{color,fancyhdr,lastpage,graphicx}
\usepackage{xspace}
\usepackage{wasysym}
\usepackage{setspace}
\usepackage{bbm}
\usepackage{stmaryrd}

\usepackage{pstricks}
\usepackage{pst-node}
\usepackage{pst-tree}

\newcommand{\NN}{\mathbb{N}}      
\newcommand{\RR}{\mathbb{R}}

   \newcommand{\DD}{\mathbb{D}}

\newcommand{\mcP}{\mathcal{P}} \newcommand{\mcC}{\mathcal{C}}  \newcommand{\mcF}{\mathcal{F}} 
   \newcommand{\mcM}{\mathcal{M}} 
        
\newcommand{\mcQ}{\mathcal{Q}}  \newcommand{\mcL}{\mathcal{L}}

  \newcommand{\be}{{\bf e}} \newcommand{\bfa}{{\bf a}}  \newcommand{\bfb}{{\bf b}} 
   
  \newcommand{\bfz}{{\bf z}}  
    \newcommand{\bfX}{{\bf X}}  \newcommand{\bfY}{\bf Y}   \newcommand{\bfB}{\bf B}

\newcommand{\EE}{\mathbb{E}} \newcommand{\PP}{\mathbb{P}}

\newcommand{\st}{such that }

\newcommand{\varep}{\varepsilon}
\newcommand{\ep}{\epsilon}
\renewcommand{\leq}{\leqslant}
\renewcommand{\geq}{\geqslant}

\newcommand{\al}{\alpha}
\newcommand{\Exp}{\textrm{exp}}

\newcommand{\wrt}{with respect to }
\renewcommand{\st}{such that }

\newcommand{\ssk}{\smallskip}

\newtheorem{thm}{\hspace{-0.15cm}  {\sc Theorem} }
\newtheorem{cor}[thm]{\hspace{-0.15cm}  {\sc Corollary} }
\newtheorem{lem}[thm]{\hspace{-0.15cm}  {\sc Lemma} }
\newtheorem{prop}[thm]{\hspace{-0.15cm} {\sc Proposition}}
\newtheorem{defn}[thm]{ \hspace{-0.3cm} {\sc Definition}}

\newtheorem{rem}[thm]{{\sc Remark}}
\newtheorem{rems}[thm]{{\sc Remarks}}

\numberwithin{equation}{section} 

\newenvironment{Dem}{%
    \begin{list}{\hspace{0.5cm}{\sc Proof --}}{%
        \setlength{\topsep}{0pt}%
        \setlength{\leftmargin}{0pt}%
        \setlength{\rightmargin}{0pt}%
        \setlength{\listparindent}{0pt}%
        \setlength{\itemindent}{0pt}%
        \setlength{\parsep}{0pt}%
        \addtolength{\leftmargin}{20pt}%
        \addtolength{\rightmargin}{0pt}%
    } \item }{\hfill{\space $\rhd$}\end{list}\smallskip}

\newenvironment{DemSewingLemma}{%
    \begin{list}{\hspace{0.5cm}{\sc Proof of theorem \ref{ThmSewingLemmaFlows} --}}{%
        \setlength{\topsep}{0pt}%
        \setlength{\leftmargin}{0pt}%
        \setlength{\rightmargin}{0pt}%
        \setlength{\listparindent}{0pt}%
        \setlength{\itemindent}{0pt}%
        \setlength{\parsep}{0pt}%
        \addtolength{\leftmargin}{20pt}%
        \addtolength{\rightmargin}{0pt}%
    } \item }{\hfill{\space $\rhd$}\end{list}\smallskip}
    
    \newenvironment{DemSewingLemmaLinear}{%
    \begin{list}{\hspace{0.5cm}{\sc Proof of theorem \ref{ThmSewingLemmaFlowsLinear} --}}{%
        \setlength{\topsep}{0pt}%
        \setlength{\leftmargin}{0pt}%
        \setlength{\rightmargin}{0pt}%
        \setlength{\listparindent}{0pt}%
        \setlength{\itemindent}{0pt}%
        \setlength{\parsep}{0pt}%
        \addtolength{\leftmargin}{20pt}%
        \addtolength{\rightmargin}{0pt}%
    } \item }{\hfill{\space $\rhd$}\end{list}\smallskip}

    \newenvironment{DemPropEstimatesGeneralRDE}{%
    \begin{list}{\hspace{0.5cm}{\sc Proof of proposition \ref{PropEstimatesGeneralRDE} --}}{%
        \setlength{\topsep}{0pt}%
        \setlength{\leftmargin}{0pt}%
        \setlength{\rightmargin}{0pt}%
        \setlength{\listparindent}{0pt}%
        \setlength{\itemindent}{0pt}%
        \setlength{\parsep}{0pt}%
        \addtolength{\leftmargin}{20pt}%
        \addtolength{\rightmargin}{0pt}%
    } \item }{\hfill{\space $\rhd$}\end{list}\smallskip}

\title{Flows driven by rough paths}
\date{\today}
\author{I. Bailleul}
\address{IRMAR, 263 Avenue du General Leclerc, 35042 RENNES, France}
\email{ismael.bailleul@univ-rennes1.fr}
\thanks{This research was partially supported by an ANR grant "Retour post-doctorant".}

\begin{document}

\maketitle

\begin{abstract}
We devise in this work a simple mechanism for constructing flows on a Banach space from approximate flows, and show how it can be used in a simple way to reprove from scratch and extend the main existence and well-posedness results for rough differential equations, in the context of dynamics on a Banach space driven by a H{\"o}lder weak geometric rough path; the explosion question under linear growth conditions, Taylor expansion and Euler estimates are also dealt with. We illustrate our approach by proving an existence and well-posedness result for some mean field stochastic rough differential equation.
\end{abstract}


\tableofcontents

\section{Introduction}

Since Lyons wrote his groundbreaking article \cite{Lyons97} on rough paths, there has been a constantly growing industry in understanding the scope of the theory. Besides providing an alternative enlighting view on Ito's theory of stochastic integration and some deep associated results, rough paths theory has now invaded the world of numerical simulations, stochastic and deterministic partial differential equations and finance, to name but a few areas. Despite Lyons' Saint Flour lecture notes \cite{LyonsStFlour}, his book \cite{LyonsQian} with Qian and the impressive and exhaustive book \cite{FVBook} of Friz and Victoir, rough path still seems to be seen as a somewhat difficult and technical subject where algebra and classical analysis meet in an intricate way.

We show in this work how the main existence and uniqueness results of the theory can be proved from scratch using a simple mechanism for constructing flows from approximate flows, in the context of dynamics in a Banach space driven by H\"older weak geometric rough paths. Contrary to Lyons, Friz-Victoir or Gubinelli's approach,  we work primarily with maps from the state space $E$ to itself rather than with $E$-valued paths. Our dynamics on the space of maps $\varphi$ will have typical form
\begin{equation}
\label{BasicRDEquation}
d\varphi = Vdt+{\textrm F}\bfX(dt),
\end{equation}
for some driving vector fields $V$ and ${\textrm F}=(V_1,\dots,V_\ell)$ and a (H\"older weak geometric $p$-) rough path $\bfX$ on some time interval $[0,T]$. Roughly speaking, a rough path consists of an $\RR^\ell$-valued non-smooth path $X$ together with a number of objects which play the role of the missing iterated integrals of $X$, in the sense that they satisfy the same algebraic relations as the iterated integrals of any smooth path, as well as some natural size requirements. Given such a rough path and some vector fields $V,V_1,\dots,V_\ell$, a family $(\mu_{ts})_{0\leq s\leq t\leq T}$ of diffeomorphisms of $E$ is constructed from an ordinary differential equation. These maps do not form a flow, in the sense that $\mu_{ts}$ is not equal to $\mu_{tu}\circ\mu_{us}$ for all $0\leq s\leq u\leq t\leq T$. However, some conditions on the driving vector fields and the rough path  ensure the existence of a unique \emph{flow} of maps $(\varphi_{ts})_{0\leq s\leq t\leq T}$ close to $(\mu_{ts})_{0\leq s\leq t\leq T}$. This is what we call in that case the solution to the rough differential equation \eqref{BasicRDEquation} on flows.

\ssk

All other approaches to rough paths consider paths associated with a point motion as the fundamental basic object. Lyons, and later Gubinelli, interpret \eqref{BasicRDEquation} as an integral equation, which requires a suitable notion of rough path integral as a mechanism to attach to a rough path $\bfY$ and some sufficiently regular one form $g$ another rough path $\int_0^{\cdot} g({\bfY})\,d{\bfY}$. Solving equation \eqref{BasicRDEquation} then amounts to find a fixed point to an integral equation of the form ${\bfY} = \int_0^\cdot g({\bfY})\,d{\bfY}$, where $\bfY$ is some extension of the original rough path $\bfX$. See for instance \cite{Lyons97}, \cite{LyonsStFlour} or \cite{LyonsQian}, and Friz and Hairer's forthcoming review of Gubinelli's approach. Davie's definition takes as a starting point the fact that solutions of a controlled ordinary differential equation $dx_t = \sum_{i=1}^\ell V_i(x_t)dh^i_t$, with $h$ smooth, admit Euler expansions
\begin{equation}
\label{EqEulerExpansion}
x_t-x_s = \sum_{i=1}^\ell (h_t-h_s)^i V_i(x-s) + \sum_{j,k=1}^\ell \Big(\int_s^t (h_r-h_s)^kdh^j_r\Big)(V_jV_k)(x_s) + o\big(|t-s|^2\big),
\end{equation}
which describe the path $x_\bullet$ accurately enough to characterize it uniquely. A $p$-rough path ${\bfX} = (X^1,X^2)$, with $2\leq p<3$, provides quantities $X^1_{ts}\in\RR^\ell$ and $X^2_{ts}\in M_\ell(\RR)$, for $0\leq s\leq t\leq T$, which, when plugged into formula \eqref{EqEulerExpansion} in place of $h_t-h_s$ and $\int_s^t (\cdots)$, with a different $o(\cdot)$ term, describe accurately a unique $\RR^d$-valued path. Friz and Victoir extended Davie's picture to any geometric  H\"older $p$-rough path by interpreting Davie's picture as a way of contructing $\RR^d$-valued paths as limits in some appropriate topology of paths generated by controlled ordinary differential equations in which the control converges in a rough path sense to some limit rough path. With such a view, no notion of integral is needed to define a dynamics. Neither is it the case in the approach developped in this work, whose core is a simple and non-commutative extension of Feyel-de la Pradelle' sewing lemma \cite{FdlP}, well-adapted to the construction of flows of maps, and totally independent of any problem about rough paths; it is the subject of section \ref{SectionConstructingFlows}. 

\ssk

It comes as a nice feature of our approach that convergence estimates for Euler/Milstein type schemes and Taylor expansion for solution flows come almost for free; this is explained in section \ref{SectionTaylorExpansion}. The approach is sufficiently robust to work with vector fields with at most linear growth and prove non-explosion of the dynamics in that case, as explained in section \ref{SectionLinearGrowth}. The results of section \ref{SectionGeneralRDE} on flows driven by rough paths are applied in section \ref{SectionPathRDE} to give a simple proof of the main existence/uniqueness results for the classical rough differential equations on paths driven by H\"older weak geometric. They are illustrated differently in section \ref{SectionNonlinearRDE} where we prove a well-posedness result for some large class of mean field stochastic rough differential equation.

\medskip

Although it is true that the classical approach to rough paths eventually leads with some work to the contruction of solution flows to equation \eqref{BasicRDEquation} in our sense, this work points out the fact that it is fruitful to take the converse direction and consider flows as basic objects rather than paths. This change of point of view has several advantages over the usual approach, the first of which being that we are able to give a simple "blackbox" for constructing flows in a Banach space, which provides an almost technical-free approach to rough differential equations. Our main results on flows driven by rough paths
\begin{itemize}
   \item well-posedness of the rough differential equation \eqref{BasicRDEquation} on flows, giving back and extending sharp existence and uniqueness results for classical rough differential equations,  
   \item non-explosion for dynamics driven by weak geometric rough paths under linear growth conditions on the vector fields,
   \item existence and well-posedness results for some mean field rough differential equations, 
\end{itemize}
hold for any weak geometric H\"older $p$-rough path on any Banach space. Our approach also offers an alternative to branched rough paths for solving rough differential equations driven by non-weak geometric rough paths and provides in that case an extension of Lyons' theorem \cite{Lyons97}, who only deals with the case $2<p<3$.

\bigskip

A few notations will be used thoughout the text, which we gather here.

\medskip

\begin{itemize}
   \item $\big(E,|\cdot|\big)$ will denote Banach space and $\textrm{L}(E)$ the set of linear continuous maps from $E$ to itself.
   \item Fix $T>0$ and a non-integer $\gamma>1$, with integer part $[\gamma]$ and fractional part $\{\gamma\}$. We say that a function or a vector field $W$ on $E$ is $\gamma$-Lipschitz if it is $\mcC^{[\gamma]}$ with a $\{\gamma\}$-H\"older continuous differential of order $[\gamma]$, and $W$ and all its derivative are bounded. Its $\gamma$-Lipschitz norm $\|W\|_\gamma$ is defined as
\begin{equation}
\label{DefGammaNorm}
\|W\|_{\gamma} = \sum_{r=0}^{[\gamma]} \big|W^{(r)}\big| +\big\|W^{([\gamma])}\big\|_{\{\gamma\}}<\infty,
\end{equation}
where $W^{(r)}$ is the differential of order $r$ of $W$, with operator norm $\big|W^{(r)}\big|$, and $\big\|W^{([\gamma])}\big\|_{\{\gamma\}}$ stands for the classical $\{\gamma\}$-H\"older norm of $W^{([\gamma])}$.   
   
   \item Given some sufficiently regular vector fields $V_1,\dots,V_\ell$ on $E$ and a tuple $I=(i_1,\dots,i_r)\in\llbracket 1,\ell\rrbracket^r$, we identify vector fields with derivation operators and write $V_I$ for the differential operator $f\in\mcC^r\mapsto V_{i_1}\big(V_{i_2}(\dots V_{i_r}f)\big)$. Writing $[V,W]$ for the bracket of two vector fields, this defines a first order differential operator, that is a vector field. Set 
   $$
   V_{[I]} = \Big[V_{i_1},\big[V_{i_2},\dots,[V_{i_{r-1}},V_{i_r}]\big]\dots\Big].
   $$
   \item We identify in the sequel $\textrm{L}\big(\RR^\ell\big)=\RR^\ell\otimes \big(\RR^\ell\big)^*$ and $\RR^\ell\otimes\RR^\ell$ with $M_\ell(\RR)$, via the matrix representation of linear maps in the canonical basis. In these terms, given any two vectors $x,y$ of $\RR^\ell$, the $(jk)$-th component $(x\otimes y)^{jk}$ of $x\otimes y$ is $x^jy^k$.
   \item We use the convention that $a^ib_i$ stands for the sum $\sum_i a^ib_i$. 
   \item Constants depending only on the appropriate norms of some given vector fields are said to depend on the data of the problem.
   \item We use the letter $c$ for constants depending on the data of the problem, and whose value may change from place to place.
   \item The notation $O_c(A)$ stands for a quantity whose norm or absolute value is bounded above by $cA$.
\end{itemize}

\medskip

\section{Constructing flows on a Banach space}
\label{SectionConstructingFlows}

Let $E$ be a Banach space. Recall that a \textit{flow} on $E$ is a family $(\varphi_{ts})_{0\leq s\leq t\leq T}$ of maps from $E$ to itself such that $\varphi_{ts}=\varphi_{tu}\,\circ\,\varphi_{us}$, for all $0\leq s\leq u\leq t\leq T$, and $\varphi_{ss}=\textrm{Id}$, for all $0\leq s\leq T$. We provide in this section a simple tool for constructing flows on $E$, which rests on an elementary extension of Feyel-de la Pradelle sewing lemma \cite{FdlP} to the non-commutative setting of maps from $E$ to itself. Given a family of maps $(\mu_{ts})_{0\leq s\leq t\leq T}$, set 
$$
\mu^n_{ts} = \bigcirc_{i=0}^{n-1}\mu_{s_{i+1}s_i},
$$ 
with $s_i=s+\frac{i}{n}\,(t-s)$. Given a partition $\pi_{ts}=\{s=s_0<s_1<\cdots<s_{n-1}<s_n=t\}$ of $(s,t)$, set 
$$
\mu_{\pi_{ts}} = \mu_{s_ns_{n-1}}\circ\cdots\circ\mu_{s_1s_0}.
$$

\begin{thm}[Sewing lemma for flows]
\label{ThmSewingLemmaFlows}
Let $\big(\mu_{ts}\big)_{0\leq s\leq t\leq T}$ be a family of Lipschitz continuous maps from $E$ to itself, depending continuously on $(s,t)$ in the uniform topology, and for which 
\begin{itemize}
   \item[{\bf H1.}] the maps $\mu_{ts} : E\rightarrow E$ have Lipschitz constant uniformly bounded above by $\big(1+o_\alpha(1)\big)$, for any $0\leq s\leq t\leq T$ with $t-s\leq \alpha$, and any $\alpha>0$. \vspace{0.1cm}
\end{itemize}
Suppose further that there exists some positive constants $c_1, c_2,\delta$ and $a>1$, such that 
\begin{itemize}
   \item[{\bf H2.}] we have
   \begin{equation}
   \label{EqMumu}
   \big\|\mu_{tu}\circ\mu_{us}-\mu_{ts}\big\|_\infty \leq c_1|t-s|^a
   \end{equation}   
   for all $0\leq s\leq u\leq t\leq T$, \vspace{0.1cm}
   \item[{\bf H3.}] the maps $\mu^n_{ts}(\cdot)$, for $n\geq 2$ and $t-s\leq \delta$, are all Lipschitz continuous, with a Lipschitz constant uniformly bounded above by $c_2$.   
\end{itemize}
Then there exists a unique flow of maps $\big(\varphi_{ts}\big)_{0\leq s\leq t\leq T}$ on $E$ such that 
   \begin{equation}
   \label{EqApproxVarphi}
   \big\|\varphi_{ts}-\mu_{ts}\big\|_\infty \leq c|t-s|^a
   \end{equation}
holds for some positive constant $c$ and all $0\leq s\leq u\leq t\leq T$, with $t-s\leq \delta$, and we have 
   \begin{equation}
   \label{EqApproxVarphiMu}
   \|\varphi_{ts}-\mu_{\pi_{ts}}\|_\infty\leq c_1c_2 T\,\big|\pi_{ts}\big|^{a-1}
   \end{equation} 
for any partition $\pi_{ts}$ of any interval $(s,t)\subset [0,T]$, of mesh $\big|\pi_{ts}\big|\leq \delta$.
\end{thm}

\ssk

A family of maps enjoying property \eqref{EqMumu} is called an {\bf approximate flow}. To prepare the proof of theorem \ref{ThmSewingLemmaFlows}, note that it is elementary to improve identity \eqref{EqMumu} under the stronger form given in lemma \ref{LemImprovedEqMumu}, and for which we introduce the following definition.

\begin{defn}
Let $\ep\in (0,1)$ be given. A partition $\pi=\{s=s_0<s_1<\cdots<s_{n-1}<s_n=t\}$ of $(s,t)$ is said to be \emph{\textbf{of special type $\ep$}} if we have $\ep\leq\frac{s_i-s_{i-1}}{s_{i+1}-s_{i-1}}\leq 1-\ep$, for all $i=1\dots n-1$. The trivial partition of any interval into the interval itself is also said to be of special type $\ep$.
\end{defn}

A partition of any interval into sub-intervals of equal length has special type $\frac{1}{2}$. Given a partition $\pi_{ts}$ of $(s,t)$ of special type $\ep$ and $u\in\{s_1,\dots,s_{n-1}\}$, the induced partitions $\pi_{us}$ and $\pi_{tu}$ of the intervals $(s,u)$ and $(u,t)$ are also of special type $\ep$. Set $c_\ep = \underset{\ep\leq\beta\leq 1-\ep}{\max}\;\big\{\beta^a+(1-\beta)^a\big\}<1$, and 
$$
L = \frac{2c_1}{1-c_\ep}.
$$

\begin{lem}
\label{LemImprovedEqMumu}
Let $\big(\mu_{ts}\big)_{0\leq s\leq t\leq T}$ be an approximate flow on $E$ satisfying assumption \emph{H1}. Given $\ep>0$, there exists a positive constant $\delta$ such that for any $0\leq s\leq t\leq T$ with $t-s\leq \delta$, and any special partition $\pi_{ts}$ of type $\ep$ of an interval $(s,t)\subset [0,T]$, we have 
\begin{equation}
\label{EqImprovedEqMumu}
\big\|\mu_{\pi_{ts}}-\mu_{ts}\big\|_\infty \leq L|t-s|^a.
\end{equation}
 \end{lem}

\begin{Dem}
We proceed by induction on the number $n$ of sub-intervals of the partition. The case $n=2$ is identity \eqref{EqMumu}. Suppose the statement has been proved for $n\geq 2$. Fix $0\leq s<t\leq T$ with $t-s\leq \delta$, and let $\pi_{ts}=\{s_0=s<s_1<\cdots<s_n<s_{n+1}=t\}$ be a partition of $(s,t)$ of special type $\ep$, splitting the interval $(s,t)$ into $(n+1)$ sub-intervals. Set $m=\big[\frac{n+1}{2}\big]$ and $u=s_m$, so the two partitions $\pi_{tu}$ and $\pi_{us}$ are both of special type $\ep$, with respective cardinals no greater than $n$, and $\ep \leq \frac{t-u}{t-s} \leq 1-\ep$. Then
\begin{equation}
\label{EqFirstEstimate}
\begin{split}
\big\|\mu_{\pi_{ts}} - \mu_{ts}\big\|_\infty &\leq \big\|\mu_{\pi_{tu}}\circ\mu_{\pi_{us}} - \mu_{tu}\circ\mu_{\pi_{us}}\big\|_\infty + \big\|\mu_{tu}\circ\mu_{\pi_{us}} - \mu_{ts}\big\|_\infty \\
&\leq \big\|\mu_{\pi_{tu}} - \mu_{tu}\big\|_\infty + \big\|\mu_{tu}\circ\mu_{\pi_{us}} - \mu_{tu}\circ\mu_{us}\big\|_\infty + \big\|\mu_{tu}\circ\mu_{us} - \mu_{ts}\big\|_\infty \\
&\leq L|t-u|^a+ \big(1+o_\delta(1)\big)L\,|u-s|^a  + c_1|t-s|^a,
\end{split}
\end{equation}
by the induction hypothesis together with assumptions H1 and H2. Set $u-s=\beta(t-s)$, with $\ep\leq \beta\leq 1-\ep$.  The above inequality rewrites 
\begin{equation*}
\begin{split}
\big\|\mu_{\pi_{ts}}-\mu_{ts}\big\|_\infty &\leq \Big\{\big(1+o_{\delta}(1)\big)\big( (1-\beta)^a+ \beta^a\big) L + c_1 \Big\} \,|t-s|^a. 
\end{split}
\end{equation*}
One closes the induction by choosing $\delta$ small enough to have $\big(1+o_\delta(1)\big)c_\ep\leq \frac{1+c_\ep}{2}$, for which choice the above term $\{\cdots\}$ is no greater than $L$.
\end{Dem}

\medskip

\begin{DemSewingLemma}
The existence and uniqueness proofs of the statement of theorem \ref{ThmSewingLemmaFlows} both rely on the elementary identity
\begin{equation}
\label{ElementaryIdentity}
f_N\circ\cdots\circ f_1 \,- \,g_N\circ\cdots\circ g_1 = \sum_{i=1}^N \Big(g_N\circ\cdots \circ g_{N-i+1}\circ f_{N-i}\, -\, g_N\circ\cdots \circ \,g_{N-i+1}\circ\, g_{N-i}\Big)\circ f_{N-i-1} \circ\cdots\circ f_1,
\end{equation}
where the $g_i$ and $f_i$ are maps from $E$ to itself, and where we use the obvious convention concerning the summand for the first and last term of the sum. In particular, if all the maps $g_N\circ\cdots\circ g_k$ are Lipschitz continuous, with a common upper bound $c'$ for their Lipschitz constants, then 
\begin{equation}
\label{CsqElementaryIdentity}
\big\|f_N\circ\cdots\circ f_1 - g_N\circ\cdots\circ g_1\big\|_\infty \leq c'\sum_{i=1}^N \|f_i-g_i\|_\infty.
\end{equation}

\ssk

\textbf{a) Existence.} Set $\textrm{D}_\delta := \big\{0\leq s\leq t\leq T\,;\,t-s\leq \delta\big\}$ and write $\DD_\delta$ for the intersection of $\textrm{D}_\delta $ with the set of dyadic real numbers. Given $s=a2^{-n_0}$ and $t=b2^{-n_0}$ in $\DD_\delta$, define for $n\geq n_0$
\begin{equation}
\label{EqMu(n)}
\mu^{(n)}_{ts} := \mu_{ts}^{2^n} = \mu_{s_{N(n)}s_{N(n)-1}}\circ\cdots\circ\mu_{s_1s_0} ,
\end{equation}
where $s_i=s+i2^{-n}$ and $s_{N(n)}=t$. Given $n\geq n_0$, write 
$$
\mu^{(n+1)}_{ts} = \overset{N(n)-1}{\underset{i=0}{\bigcirc}}\big(\mu_{s_{i+1}s_i+2^{-n-1}}\circ\mu_{s_i+2^{-n-1}s_i}\big)
$$
and use \eqref{ElementaryIdentity} with $f_i = \mu_{s_{i+1}s_i+2^{-n-1}}\circ\mu_{s_i+2^{-n-1}s_i}$ and $g_i=\mu_{s_{i+1}s_i}$ and the fact that all the maps $\mu_{s_{N(n)}s_{N(n)-1}} \circ \cdots \circ \mu_{s_{N(n)-i+1} s_{N(n)-i}} = \mu_{s_{N(n)}s_{N(n)-i}}^i$ are Lipschitz continuous with a common Lipschitz constant $c_2$, by assumption H3, to get by \eqref{CsqElementaryIdentity} and \eqref{EqMumu}
$$
\Big\|\mu^{(n+1)}_{ts} - \mu^{(n)}_{ts}\Big\|_\infty \leq c_2\sum_{i=0}^{N(n)-1}\big\|\mu_{s_{i+1}s_i+2^{-n-1}}\circ\mu_{s_i+2^{-n-1}s_i} - \mu_{s_{i+1}s_i}\big\|_\infty \leq  c_1c_2 T\,2^{-(a-1)n};
$$
so $\mu^{(n)}$ converges uniformly on $\DD_\delta$ to some continuous function $\varphi$. We see that $\varphi$ satisfies inequality \eqref{EqApproxVarphi} on $\DD_\delta$ as a consequence of \eqref{EqImprovedEqMumu}. As $\varphi$ is a uniformly continuous function of $(s,t)\in\DD_\delta$, by \eqref{EqApproxVarphi}, it has a unique continuous extension to $\textrm{D}_\delta$, still denoted by $\varphi$. To see that it defines a flow on $\textrm{D}_\delta$, notice that for dyadic times $s\leq u\leq t$, we have $\mu^{(n)}_{ts} = \mu^{(n)}_{tu}\circ\mu^{(n)}_{us}$, for $n$ big enough; so $\varphi_{ts}=\varphi_{tu}\circ\varphi_{us}$ for such triples of times in $\DD_\delta$, hence for all times since $\varphi$ is continuous. The map $\varphi$ is easily extended as a flow to the whole of $\big\{(s,t)\in [0,T]^2,\;\,0\leq s\leq t\leq T\big\}$.

\medskip

\noindent \textbf{b) Uniqueness.} Let $\psi$ be any flow satisfying condition \eqref{EqApproxVarphi}. With formulas \eqref{ElementaryIdentity} and \eqref{CsqElementaryIdentity} in mind, rewrite \eqref{EqApproxVarphi} under the form $\psi_{ts}=\mu_{ts}+O_c\bigl(|t-s|^a\bigr)$, with obvious notations. Then
\begin{equation*}
\begin{split}
\psi_{ts} &= \psi_{s_{2^n}s_{2^n-1}}\circ\cdots\circ\psi_{s_1s_0} = \Bigl(\mu_{s_{2^n}s_{2^n-1}} + O_c\bigl(2^{-an}\bigr)\Bigr) \circ\cdots\circ \Bigl(\mu_{s_1s_0} + O_c\bigl(2^{-an}\bigr)\Bigr)  \\
             &= \mu_{s_{2^n}s_{2^n-1}}\circ\cdots\circ\mu_{s_1s_0} + \Delta_n = \mu^{(n)}_{ts} + \Delta_n,
\end{split}
\end{equation*}
where $\Delta_n$ is of the form of the right hand side of \eqref{ElementaryIdentity}, so is bounded above by a constant multiple of $2^{-(a-1)n}$, since all the  maps $\mu_{s_{2^n}s_{2^n-1}} \circ \cdots \circ \mu_{s_{2^n-\ell+1} s_{2^n-\ell}}$ are Lipschitz continuous with a common upper bound for their Lipschitz constants, by assumption H3.  Sending $n$ to infinity shows that $\psi_{ts}=\varphi_{ts}$.

\medskip

\noindent {\bf c) Speed of convergence.} Given any partition $\pi_{ts} = \{s_0=s<\cdots<s_n=t\}$ of $(s,t)$, using as above the uniform Lipschitz character of  the maps $\mu_{s_ns_{n-1}} \circ \cdots \circ \mu_{s_{i+s} s_i}$, and writing $\varphi_{ts} = \bigcirc_{i=0}^{n-1} \varphi_{s_{i+1}s_i}$, we see as a consequence of \eqref{CsqElementaryIdentity} that we have for $\big|\pi_{ts}\big|\leq \delta$
\begin{equation*}
\big\|\varphi_{ts}-\mu_{\pi_{ts}}\big\|_\infty \leq c_2\sum_{i=0}^{n-1} \big\|\varphi_{s_{i+1}s_i}-\mu_{s_{i+1}s_i}\big\|_\infty \leq c_1c_2\sum_{i=0}^{n-1} |s_{i+1}-s_i|^a \leq c_1c_2T\,\big|\pi_{ts}\big|^{a-1}.
\end{equation*}
\end{DemSewingLemma}

The next proposition provides a simple mean for obtaining the uniform control of the Lipschitz size of the maps $\mu^n_{ts}$, needed to apply the sewing lemma for flows; it is a simple variation of lemma \ref{LemImprovedEqMumu}. We write $|M|$ for the operator norm of a linear map $M : E\rightarrow E$.

\begin{prop}[Uniform Lipschitz controls]
\label{PropUniformLipControl}
Let $\alpha$ and $\rho$ be positive constants, with $0<1-\rho<\alpha<1$, and 
\begin{itemize}
   \item[{\bf H0.}] $(\mu_{ts})_{0\leq s\leq t\leq T}$ be an approximate flow of $(1+\rho)$-Lipschitz maps from $E$ to itself such that one can write
   \begin{equation}
   \label{EqH0}
   \quad D_x\mu_{ts} = \textrm{\emph{Id}} + A^{ts}_x + B^{ts}_x,
   \end{equation} 
\end{itemize}
for some $\textrm{\emph{L}}(E)$-valued $\rho$-Lipschitz maps $A^{ts}$ with $\rho$-Lipschitz norm bounded above by $c|t-s|^\alpha$, and some $\textrm{\emph{L}}(E)$-valued $\mcC^1$ bounded maps $B^{ts}$, with $\mcC^1$-norm bounded above by $o_{t-s}(1)$. If there exists a positive constants $c_3$ such that we have for all $x\in E$
\begin{equation}
\label{EqUniformLipControls}
\big|D_x(\mu_{tu}\circ\mu_{us}) - D_x\mu_{ts}\big| \leq c_3 |t-s|^a
\end{equation}
for all $0\leq s\leq t\leq T$, with the same constant $a>1$ as in \eqref{EqMumu}, then, given $\ep>0$, there exists two positive constants $\delta$ and $K$ such that 
$$
\big|D_x(\mu_{\pi_{ts}}) - D_x\mu_{ts}\big| \leq K |t-s|^a
$$
holds for any partition $\pi_{ts}$ of $(s,t)$ of special type $\ep$ and all $x\in E$, whenever $t-s\leq \delta$.
\end{prop}

\begin{Dem}
We proceed by induction on the number $n$ of sub-intervals of the partition as in the proof of lemma \ref{LemImprovedEqMumu}. The case $n=2$ is identity \eqref{EqUniformLipControls}. Suppose the statement has been proved for $n\geq 2$. Fix $0\leq s<t\leq T$ with $t-s\leq \delta$, and let $\pi_{ts}=\{s_0=s<s_1<\cdots<s_n<s_{n+1}=t\}$ be a partition of $(s,t)$ of special type $\ep$, splitting the interval $(s,t)$ into $(n+1)$ sub-intervals. Set $m=\Big[\frac{n+1}{2}\Big]$ and $u:=s_m$, so the two partitions $\pi_{tu}$ and $\pi_{us}$ are both of special type $\ep$, with respective cardinals no greater than $n$. Then we have for any $x\in E$
\begin{equation*}
\label{FirstEstimate}
\begin{split}
D_x\big(&\mu_{\pi_{ts}}\big) - D_x\mu_{ts} = D_x\big(\mu_{\pi_{tu}}\circ\mu_{\pi_{us}}\big) - D_x\mu_{ts} \\
        													   &= \Big(D_{\mu_{\pi_{us}}(x)}\mu_{\pi_{tu}} - D_{\mu_{\pi_{us}}(x)}\mu_{tu} \Big)\big(D_x\mu_{\pi_{us}}\big) + \big(D_{\mu_{\pi_{us}}(x)}\mu_{tu} - D_{\mu_{us}(x)}\mu_{tu}\big)\big(D_x\mu_{\pi_{us}}\big) \\
        													   &+\big(D_{\mu_{us}(x)}\mu_{tu}\big)\Big(D_x\mu_{\pi_{us}}-D_x\mu_{us}  \Big) + \Big(\big(D_{\mu_{us}(x)}\mu_{tu}\big)\big(D_x\mu_{us}\big) - D_x\mu_{ts}\Big) \\
        													   &=: (1) + (2) + (3) + (4).
\end{split}
\end{equation*}
We treat each term separately using repeatedly the induction hypothesis and lemma \ref{LemImprovedEqMumu} when needed. We first have 
$$
\big|(1)\big| \leq \big(1+o_\delta(1)+\delta^a K\big) K|t-u|^a.
$$ 
Also, using the fact that $|t-u|\leq \frac{1-\ep}{\ep}|u-s|$, we see by \eqref{EqH0} that we have
$$
\Big|D_{\mu_{\pi_{us}}(x)}\mu_{tu} - D_{\mu_{us}(x)}\mu_{tu}\Big| \leq c|t-u|^\alpha L^\rho|u-s|^{a\rho} + o_{|t-u|}(1)L|u-s|^a \leq o_\delta(1)|u-s|^a,
$$
provided $a<\alpha+\rho a$, which we can suppose without loss of generality since $1-\rho<\alpha<1$. As the term $D_x\mu^m_{us}$ has size no greater than $\big(1+o_\delta(1)\big) + K|u-s|^a$, we have
$$
\big|(2)\big| \leq \big\{\big(1+o_\delta(1)\big) + K|u-s|^a\big\}\,o_\delta(1)|u-s|^a =: c_Ko_\delta(1)|u-s|^a.
$$
Last, we have the upper bound
$$
\big|(3)\big| \leq \big(1+o_\delta(1)\big) K|u-s|^a,
$$
while  $\big|(4)\big| \leq c_3|t-s|^a$, by \eqref{EqUniformLipControls}. All together, and writing $t-u=\beta (t-s)$, with $\ep\leq beta\leq 1-\ep$, this gives 
\begin{equation*}
\begin{split}
\big|D_x\big(\mu_{\pi_{ts}}\big) - D_x\mu_{ts}\big| &\leq \Big(\Big\{\big(1+o_\delta(1)+\delta^a K\big)\beta^a + \big(1+o_\delta(1)\big)(1-\beta)^a\Big\}K + c_3+c_Ko_\delta(1)\Big)|t-s|^a \\
&\leq \Big\{\big(1+o_\delta(1)+\delta^a K\big)c_\ep K + \big(c_3+c_Ko_\delta(1)\big)\Big\}|t-s|^a
\end{split}
\end{equation*}
The induction is closed by choosing $\delta$ and $K$ so as to have $(\cdots)c_\ep<1$ and $\big\{(\cdots)c_\ep K + (\cdots)\big\}\leq K$, in the above upper bound.
\end{Dem}

\ssk

An approximate flow satisfying the regularity assumption H0 and inequality \eqref{EqUniformLipControls} is said to be a {\bf $\mcC^1$-approximate flow}. We obtain our main tool for constructing flows on $E$ as a consequence of the sewing lemma for flows and proposition \ref{PropUniformLipControl}.

\begin{thm}[Contructing flows on $E$]
\label{CorMainWorkingResult}
A $\mcC^1$-approximate flow defines a unique flow $(\varphi_{ts})_{0\leq s\leq t\leq T}$ on $E$ such that 
\begin{equation}
\label{EqConstructingFlows}
\big\|\varphi_{ts}-\mu_{ts}\big\|_\infty \leq c|t-s|^a
\end{equation}
holds for some $c$, for all $0\leq s\leq t\leq T$ sufficiently close; this flow satisfies the approximation inequality \eqref{EqApproxVarphiMu}.
\end{thm}

\ssk

As a straightforward application, it is elementary to see that one defines a $\mcC^1$ approximate flow setting 
\begin{equation}
\label{EqDefnMutsODE}
\mu_{ts}(x) = x + \big(h_t-h_s\big)^iV_i(x),
\end{equation}
where $h : [0,T]\rightarrow \RR^\ell$ is a $\mcC^1$ control and the vector fields $V_1,\dots,V_\ell$ on $E$ are $\mcC^2_b$. Inequalities \eqref{EqMumu} and \eqref{EqUniformLipControls} hold in that case with $a=2$. The associated flow coincides with the flow associated with the controlled ordinary differential equation 
\begin{equation}
\label{EqControlledODE}
dx_t = V_i(x_t)dh^i_t,
\end{equation}
as the latter obviously satisfies inequality \eqref{EqConstructingFlows}.

\section{Flows driven by H\"older weak geometric $p$-rough paths}
\label{SectionGeneralRDE}

H\"older $p$-rough paths, which control the rough differential equation \eqref{BasicRDEquation} and play the role of $h$ in \eqref{EqControlledODE}, are defined in section \ref{SectionAlgebraicInterlude}. As $\RR^\ell$-valued paths, they are not regular enough for formula \eqref{EqDefnMutsODE} to define an approximate flow. The missing bit of information needed to stabilize the situation is a substitute of the non-existing iterated integrals $\int_s^t X^j_rdX^k_r$, and higher order iterated integrals, which provide a partial description of what happened to $X$ during any time interval $(s,t)$. The higher order parts of a $p$-rough path provide precisely that information.  It is an important fact that $p$-rough paths take values in a very special kind of algebraic structure, of which we recall the basic features in section \ref{SectionAlgebraicInterlude}. We shall then see in section \ref{SectionGeneralCase} how to associate to a rough path and some smooth enough vector fields a $\mcC^1$-approximate flow.

\subsection{An algebraic prelude: tensor algebra over $\RR^\ell$ and free nilpotent Lie group}
\label{SectionAlgebraicInterlude}

For $N\in\NN\cup\{\infty\}$, write $T^{(N)}_\ell$ for the direct sum $\underset{r=0}{\overset{N}{\bigoplus}}\big(\RR^\ell\big)^{\otimes r}$, with the convention that $\big(\RR^\ell\big)^{\otimes 0}$ stands for $\RR$. Denote by $\bfa =  \underset{r=0}{\overset{N}{\oplus}} a_r$ and $\bfb =  \underset{r=0}{\overset{N}{\oplus}} b_r$ two generic elements of $T^{(N)}_\ell$. The vector space $T^{(N)}_\ell$ is an algebra for the operations
\begin{equation}
\label{AlgebraOperations}
\begin{split}
&\bfa + \bfb = \underset{r=0}{\overset{N}{\oplus}} (a_r+b_r), \\
&\bfa\bfb = \underset{r=0}{\overset{N}{\oplus}} c_r, \quad\textrm{with }\; c_r=\sum_{k=0}^r a_k\otimes b_{r-k}
\end{split}
\end{equation}
It is called the \textbf{(truncated) tensor algebra of $\RR^\ell$} (if $N$ is finite). 

\ssk

The exponential map $\exp : T^{(\infty)}_\ell\rightarrow T^{(\infty)}_\ell$ and the logarithm map $\log : T^{(\infty)}_\ell\rightarrow T^{(\infty)}_\ell$ are defined by the usual series
\begin{equation}
\label{EqExpLog}
\exp(\bfa) = \sum_{n\geq 0}\frac{\bfa^n}{n !},\quad \log(\bfb) = \sum_{n\geq 1}\frac{(-1)^n}{n}(1-\bfb)^n,
\end{equation}
with the convention $a^0=1\in \RR\subset T^{(\infty)}_\ell$. Denote by $\pi_N : T^{(\infty)}_\ell\rightarrow T^{(N)}_\ell$ the natural projection. We also denote by $\exp$ and $\log$ the restrictions to $T^{(N)}_\ell$ of the maps $\pi_N\circ\exp$ and $\pi_N\circ\log$ respectively. Denote by $T^{(N),1}_\ell$, resp. $T^{(N),0}_\ell$, the elements $a_0\oplus\cdots\oplus c_N$ of $T^{(N)}_\ell$ \st $a_0=0$, resp. $a_0=1$. All the elements of $T^{(N),1}_\ell$ are invertible, and $\exp : T^{(N),0}_\ell\rightarrow T^{(N),1}_\ell$ and $\log : T^{(N),1}_\ell\rightarrow T^{(N),0}_\ell$ are reciprocal bijections.

\ssk

The formula $[\bfa,\bfb]=\bfa\bfb-\bfb\bfa$, defines a Lie bracket on $T^{(N)}_\ell$. Define inductively $F=F^1 = \RR^\ell$, considered as a subset of $T^{(\infty)}_\ell$, and $F^{n+1}=[F,F^n]\subset T^{(\infty)}_\ell$.

\begin{defn}
\begin{itemize}
   \item The Lie algebra $\frak{g}^{N}_\ell$ generated by the $F^1,\dots,F^N$ is called the \emph{\textbf{$N$-step free nilpotent Lie algebra}}.
   \item As a consequence of Baker-Campbell-Hausdorf-Dynkin formula, the subset $\exp\big(\frak{g}^{N}_\ell\big)$ of $T^{(N)}_\ell$ is a group for the multiplication operation. It is called the \emph{\textbf{$N$-step nilpotent Lie group}} on $\RR^\ell$ and denoted by $G^{(N)}_\ell$.
\end{itemize}
\end{defn}

Note that the restriction to $G^{(q)}_\ell$ of the projection map $\pi_{pq} : T^{(q)}_\ell\rightarrow T^{(p)}_\ell$, sending $a_0\oplus\cdots\oplus a_q$ to $a_0\oplus\cdots\oplus a_p$, provides a natural projection $\pi_{pq}$ from $G^{(q)}_\ell$ to $G^{(p)}_\ell$ for any $p<q$.

\medskip

The relevance of the algebraic framework provided by the $N$-step nilpotent Lie group for the study of smooth paths was first noted by Chen in his seminal work \cite{Chen77}. Indeed, for any $\RR^\ell$-valued smooth path $(x_s)_{s\geq 0}$, the family of iterated integrals 
$$
{\frak X}^N_{ts} := \left(1,x_t-x_s,\int_s^t\int_s^{s_1}dx_{s_2}\otimes dx_{s_1},\dots,\int_{s\leq s_1\leq \cdots\leq s_N\leq t}dx_{s_1}\otimes\cdots\otimes dx_{s_N}\right)
$$
defines for all $0\leq s\leq t$ an element of $G^{(N)}_\ell$. It suffices to notice that, as a function of $t$, the function ${\frak X}^N_{ts}$ satisfies the differential equation 
$$
d{\frak X}^N_{ts} = {\frak X}^N_{ts}\otimes dx_t,
$$
in $T^{(N)}_\ell$. Rough paths and weak geometric rough paths are somehow an abstract version of these objects.

\begin{defn}
Let $2\leq p$. A \emph{\textbf{H\"older $p$-rough path on $[0,T]$}} is a $T_\ell^{([p]),1}$-valued path ${\bfX} : t\in [0,T]\mapsto 1\oplus X^1_t\oplus X^2_t\oplus\cdots\oplus X^{[p]}_t$ \st 
\begin{equation}
\label{ConditionsHolder}
\underset{0\leq s<t\leq T}{\sup}\,\frac{|X^i_{ts}|}{|t-s|^{\frac{i}{p}}}<\infty,
 \end{equation}
for all $i=1\dots [p]$, where we set ${\bf X}_{ts} = {\bf X}_s^{-1}{\bf X}_t$. We define the norm of $\bfX$ to be 
\begin{equation}
\|{\bfX}\| = \underset{i=1\dots [p]}{\max}\;\underset{0\leq s<t\leq T}{\sup}\,\frac{|X^i_{ts}|}{|t-s|^{\frac{i}{p}}},
\end{equation}
and a distance $d(\bf X,\bf Y) = \|\bfX-\bfY\|$ on the set of H\"older $p$-rough path. A \emph{\textbf{H\"older weak geometric $p$-rough path on $[0,T]$}} is a $G_\ell^{([p])}$-valued $p$-rough path.
\end{defn}

\ssk

\noindent For $2\leq p<3$, the relation ${\bfX}_{us}{\bfX}_{tu}={\bfX}_{ts}$, for $0\leq s\leq u\leq t\leq T$, is equivalent to 
\begin{itemize}
   \item[{\bf (i)}] $X^1_{ts} = X^1_{tu}+X^1_{us}$,
   \item[{\bf (ii)}]  $X^2_{ts} = X^2_{tu} + X^1_{us}\otimes X^1_{tu} + X^2_{us}$.
\end{itemize}
Condition (i) means that $X^1_{ts}=X^1_{t0}-X^1_{s0}$ represents the increment of the $\RR^d$-valued path $\big(X^1_{r0}\big)_{0\leq r\leq T}$. Condition (ii) is nothing but the analogue of the elementary property $\int_s^t\int_s^r =\int_s^u\int_s^r + \int_u^t\int_s^u + \int_u^t\int_u^r$, satisfied by any reasonable notion of integral on $\RR$ such that $\int_s^t = \int_s^u+\int_u^t$. This remark justifies thinking of the $\big(\RR^\ell\otimes\RR^\ell\big)$-part of a rough path as a kind of iterated integral of $X^1$ against itself. In that setting, a $p$-rough path $\bfX$ is a weak geometric $p$-rough path iff the symmetric part of $X^2_{ts}$ is $\frac{1}{2}X^1_{ts}\otimes X^1_{ts}$, for all $0\leq s\leq t\leq T$.

\medskip

\begin{rem}
H\"older $p$-rough paths appear naturally in a probabilistic context. Let $B$ be a Brownian motion; the random process ${\bfB} = \Big(B_t-B_s,\,\int_s^tB_r\otimes {\circ d}B_r\Big)_{0\leq s\leq t\leq T}$ is almost-surely a weak geometric H\"older $p$-rough path, for any $2< p<3$. It is called the {\bf Brownian rough path}. Note that using an Ito integral, the formula ${\bf B}^\textrm{I} = \Big(B_t-B_s,\,\int_s^tB_r\otimes dB_r\Big)_{0\leq s\leq t\leq T}$, defines a H\"older $p$-rough path which is not weak geometric.
\end{rem}

\ssk

Denote by $(e_1,\dots,e_\ell)$ the canonical basis of $\RR^\ell\subset T^{(\infty)}_\ell$ and write for a tuple $I=(i_1,\dots,i_r)$
$$
{\be}_{[I]}=\Big[e_{i_1},\big[e_{i_2},\dots[e_{i_{r-1}},e_{i_r}]\big]\dots\Big],\quad \textrm{ and   } {\be}_I=e_{i_1}e_{i_2}\cdots e_{i_r},
$$
where the above products are in $T^{(\infty)}_\ell$. Write $Y^{r,I}$, with $|I|=r$, or simply $Y^I$, for the coordinates of an element $Y$ of $T_\ell^{[p]}$ in its canonical basis.

\ssk

Given a H\"older weak geometric $p$-rough path $\bfX$, denote by ${\bf\Lambda} = 0\oplus\Lambda^1\oplus\cdots\oplus\Lambda^{[p]}$ its logarithm, in Magnus-Chen-Strichartz form \cite{Strichartz}, 
\begin{equation}
\label{LogarithmIdentity}
\exp\left(\Lambda^I_{ts}{\be}_{[I]} \right) =  X^I_{ts}{\be}_I = {\bfX}_{ts},
\end{equation}
for all $0\leq s\leq t\leq T$; it takes values in the finite dimensional Lie algebra $\frak{g}^{[p]}_\ell$. Notice that since ${\bf \Lambda}$ is polynomial in $\bfX$, by formula \eqref{EqExpLog}, it is a continuous function of $\bfX$.

\subsection{Flows driven by H\"older weak geometric $p$-rough paths}
\label{SectionGeneralCase}

Given a bounded Lipschitz continuous vector field $V$, and some $\gamma$-Lipschitz vector fields $V_1,\dots,V_\ell$ on $E$, let $\mu_{ts}$ be the well-defined time $1$ map associated with the ordinary differential equation 
\begin{equation}
\label{DeterministicODEGeneralCase}
\dot y_u = (t-s)V(y_u) + \sum_{r=1}^{[p]}\sum_{I\in\llbracket 1,\ell\rrbracket^r} \Lambda_{ts}^{r,I} V_{[I]}(y_u), \quad 0\leq u\leq 1.
 \end{equation}
This equation is to be understood as the classical ordinary differential equation version of equation \eqref{BasicRDEquation}; the definition of a solution flow to equation \eqref{BasicRDEquation} given below will make that point clear. The property of $\mu_{ts}$  proved in proposition \ref{PropEstimatesGeneralRDE} below is the main reason for its introduction; it roughly says that $\mu_{ts}$ has the awaited Euler expansion, in accordance with what happens for ordinary differential equations, as emphasized in equation \eqref{EqEulerExpansion}. For $2\leq p<3$, and ${\bf X} = (X,\mathbb{X})$, equation \eqref{DeterministicODEGeneralCase} reads
$$
\dot y_u = (t-s)V(y_u) + X^i_{ts}V_i(y_u) + \frac{1}{2}\Big\{\mathbb{X}_{ts}+\frac{1}{2}X_{ts}\otimes X_{ts}\Big\}^{jk}[V_j,V_k](y_u).
$$
As the matrix $X_{ts}\otimes X_{ts}$ is symmetric, we have $X^j_{ts}X^k_{ts}\,[V_j,V_k]=0$, so \eqref{DeterministicODEGeneralCase} simplifies into 
$$
\dot y_u = (t-s)V(y_u) + X^i_{ts}V_i(y_u) + \frac{1}{2}\mathbb{X}_{ts}^{jk}[V_j,V_k](y_u).
$$

\medskip

\noindent If the vector field $V$ is $\mcC^1$ and the $V_1,\dots,V_\ell$ are $\gamma$-Lipschitz, the classical results on the dependence of solutions to ordinary differential equation \wrt some parameters ensure that for any reals $a,b^I$, the map $\Exp\bigl(aV+b^IV_{[I]}\bigr)(x)$ associated with the differential equation 
$$
\frac{d}{du}\,y_u = aV(y_u) + b^IV_{[I]}(y_u),\quad 0\leq u\leq 1, 
$$
is a continuously differentiable function of $(a,b)$. The following basic fact comes as a consequence of the analytic properties of any H\"older $p$-rough path and the definition of the topology on the set of H\"older $p$-rough paths. We write $\mu_{ts}^{\bf X}$ instead of $\mu_{ts}$ in the proposition below to emphasize its dependence on $\bf X$.

\begin{prop}
\label{PropContinuityFlow}
\begin{enumerate}
   \item Suppose $V$ is $\mcC^1$ and the $V_i$ are $\gamma$-Lipschitz. Then the diffeomorphisms $\big(\mu_{ts}^{\bf X}\big)_{0\leq s\leq t\leq T}$ depend continuously on $\big((s,t),\bfX\big)$, in the sense that 
$$
\max\Big\{\|\mu_{ts}^{\bf X}-\mu_{vu}^{\bf Y}\|_\infty,\big\|\big(\mu_{ts}^{\bf X}\big)^{-1}-\big(\mu_{vu}^{\bf Y}\big)^{-1}\big\|_\infty\Big\} \leq f\big((s,t),(u,v)\,;\,{\bf X},{\bf Y}\big)
$$
for a continuous function $f$ with $f\big((s,t),(s,t)\,;\,{\bf X,X}\big) = 0$, for all $0\leq s\leq t\leq T$, and any H{\"o}lder $p$-rough path $\bf X$. \vspace{0.2cm}

   \item Suppose $V$ is Lipschitz continuous and the $V_i$ are $(\gamma-1)$-Lipschitz. Then all the $\mu_{ts}$ are uniformly Lipschitz continuous, with Lipschitz constant no greater than $1+c|t-s|^{1/p}$, for a constant $c=c'\,\big(1+\|\bfX\|^\gamma\big)$, with $c'$ depending onthe data of the problem.
\end{enumerate}
\end{prop}

\medskip

For $(\gamma-1)$-Lipschitz vector fields $V_i$, the differential equation \eqref{DeterministicODEGeneralCase} is not necessarily well-posed. Peano's theorem ensures however the existence of a solution in a finite dimensional setting; $\mu_{ts}(x)$ stands in that case for the time $1$ value of an arbitrarily chosen solution started from $x$. The situation is different in an infinite dimensional setting where Peano's theorem may not hold without further conditions on the driving vector fields \cite{Bogachev}. There are two classical ways to get sorted out, by either imposing condition on the Banach space, or by requiring some compactness properties from the vector fields. Bownds and Diaz \cite{BowndsDiaz} prove for instance that if a continuous vector field $W : E\rightarrow E$ maps a ball $\{x\in E\,;\,|x-x_0|\leq R\}$ into a compact set contained in $\{y\in E\,\;\,|y|\leq M\}$, then the differential equation $x'_t=W(x_t)$, started from $x_0$, has a solution defined on the time interval $\big(-\frac{R}{M},\frac{R}{M}\big)$. The compactness assumption on the range of $W$ can be relaxed provided there is a good control of a measure of its weak non-compactness \cite{CramerLakshMitchell}. To fix a precise setting, we put forward the following condition.  \vspace{0.1cm} \\

\noindent {\bf (C)} \textit{The vector fields $V_1,\dots,V_\ell$ and $V_{[I]}$, with $I\in\llbracket 1,\ell\rrbracket^r$ and $r\leq [\gamma],$ map any ball $\big\{y\in E\,;\,|y-y_0|\leq R\big\}$ into a compact set.} \\ 

Assumption {\bf (C)} automatically holds in a finite dimensional setting. A basic scaling argument justifies that equation \eqref{DeterministicODEGeneralCase} has a solution defined up to time 1 provided $(t-s)$ is sufficiently small. Given an initial point $x\in E$, let $(y_u)_{0\leq u\leq 1}$ stand for any of them, and set $\mu_{ts}(x) = y_1$. 

\medskip

The following proposition is our basic step for studying flows driven by rough paths. 

\begin{prop}
\label{PropEstimatesGeneralRDE}
Let $0<\rho<1$ and $2\leq p<\gamma< [p]+1$ be given. Let $V$ be a Lipschitz-continuous vector field on $E$, and $V_1,\dots,V_\ell$ be $(\gamma-1)$-Lipschitz vector fields on $E$ satisfying the compactness assumption {\bf (C)}. Let $\bfX=1\oplus X^1\oplus\cdots\oplus X^{[p]}$ be a H\"older weak geometric $p$-rough path. Then there exists a positive constant $c$, depending only on the data of the problem \st 
\begin{equation}
\label{GralFundamentalEstimate}
 \Big\|f\circ\mu_{ts} - \Bigl\{f + (t-s)Vf + \sum_{r=1}^{[p]} \sum_{I\in\llbracket 1,\ell\rrbracket^r} X^{r,I}_{ts}V_If\Bigr\}\Big\|_\infty \leq c\big(1+\|\bfX\|^{\gamma}\big)\,\|f\|_\gamma\,|t-s|^{\frac{\gamma}{p}}
\end{equation}
holds for any $f\in\mcC^\gamma$. 
\end{prop}

The proof of this proposition and the following one are based on the following elementary identity, obtained by applying repeatedly the identity
$$
f(y_r) = f(x) + (t-s)\int_0^1(Vf)(y_u)\,du + \sum_{I}\Lambda^I_{ts}\int_0^1\big(V_{[I]}f\big)(y_u)\,du,
$$
and by separating the terms according to their size in $|t-s|$. Set $\Delta_n := \big\{(s_1,\dots,s_n)\in [0,T]^n\,;\,s_1\leq \cdots\leq s_n\big\}$, for $2\leq n\leq [p]$. For a $\gamma$-Lipschitz function $f$,  we have
\begin{equation}
\label{EqExactFormulaMuTs}
\begin{split}
f\big(\mu_{ts}(x)\big) &= f(x) + (t-s)\big(Vf\big)(x) + \sum_{k=1}^n \frac{1}{k!}\sum_{|I_1|+\dots+|I_k|\leq [p]} \left(\prod_{m=1}^k \Lambda^{I_m}_{ts}\right) \big(V_{[I_k]}\cdots V_{[I_1]}f\big)(x) \\
&+\sum_{|I_1|+\dots+|I_n|\leq [p]} \left(\prod_{m=1}^n \Lambda^{I_m}_{ts}\right) \int_0^1 \Big\{\big(V_{[I_n]}\cdots V_{[I_1]}f\big)(y_{s_n})-\big(V_{[I_n]}\cdots V_{[I_1]}f\big)(x)\Big\}{\bf 1}_{\Delta_n}\,ds_n\dots ds_1 \\
&+ (t-s)\int_0^1 \big\{\big(Vf\big)(y_r)-\big(Vf\big)(x)\big\}dr \\
&+ \sum_{k=1}^{n-1}\sum_{|I_1|+\dots+|I_k|\geq [p]+1} \left(\prod_{m=1}^k \Lambda^{I_m}_{ts}\right) \big(V_{[I_k]}\cdots V_{[I_1]}f\big)(x) \\
&+ (t-s)\sum_{k=1}^{n-1}\sum_{I_1,\dots,I_k} \left(\prod_{m=1}^k \Lambda^{I_m}_{ts}\right) \int_0^1 \big(VV_{[I_k]}\cdots V_{[I_1]}f\big)(y_{s_k}){\bf 1}_{\Delta_k}\,ds_k\dots ds_1 \\
&+ \sum_{|I_1|+\dots+|I_n|\geq [p]+1} \left(\prod_{m=1}^k \Lambda^{I_m}_{ts}\right) \int_0^1 \Big\{\big(V_{[I_n]}\cdots V_{[I_1]}f\big)(y_{s_n})-\big(V_{[I_n]}\cdots V_{[I_1]}f\big)(x)\Big\}{\bf 1}_{\Delta_n}\,ds_n\dots ds_1
\end{split}
\end{equation}
We denote by $\ep^{f\,;\,n}_{ts}(x)$ the sum of the last four lines, made up of terms of size at least $|t-s|^\frac{\gamma}{p}$. Note that this formula makes sense for all $2\leq n\leq [p]$, for $(\gamma-1)$-Lipschitz vector fields $V_i$. In the case where $n=[p]$, the terms in the second line involve only indices $I_j$ with $|I_j|=1$, so the elementary estimate 
\begin{equation}
\label{EqElementaryEstimateYu}
\big|y_r-x\big| \leq c\big(1+\|{\bfX}\|^\gamma\big)|t-s|^{1/p}, \quad 0\leq r\leq 1,
\end{equation}
can be used to control the increment in the integral, showing that this second line is of order $|t-s|^\frac{\gamma}{p}$, as the maps $V_{[I_{[p]}]}\dots V_{[I_1]}f$ are $(\gamma-[p])$-Lipschitz; we include it in the remainder $\ep^{f\,;\,[p]}_{ts}(x)$.

\medskip

\begin{DemPropEstimatesGeneralRDE}
Applying the above formula for $n=[p]$, together with the fact that $\exp({\bf \Lambda}) = \bfX$, we get the identity
\begin{equation*}
f\big(\mu_{ts}(x)\big) = f(x) + (t-s)\big(Vf\big)(x) + \sum_{I} X^I_{ts}\big(V_If\big)(x) + \ep^{f\,;\,[p]}_{ts}(x).
\end{equation*}
It is clear on the formula for $\ep^{f\,;\,[p]}_{ts}(x)$ that its absolute value is bounded above by a constant multiple of $\big(1+\|{\bfX}\|^\gamma\big)|t-s|^\frac{\gamma}{p}$, for a constant depending only on the data of the problem and $f$ as in \eqref{GralFundamentalEstimate}. 
\end{DemPropEstimatesGeneralRDE}

\medskip

A further look at formula \eqref{EqExactFormulaMuTs} also makes it clear that if $V$ is $(1+\rho)$-Lipschitz and the $V_i$ are $\gamma$-Lipschitz, with $f$ sufficiently regular, then the remainders $\ep^{f\,;\,n}_{ts}, 2\leq n\leq [p]$, have the norm of their first derivative bounded above by a constant multiple of $\big(1+\|{\bfX}\|^\gamma\big)|t-s|^a$, for a constant depending only on the data of the problem, and $a=\min\Big\{1+\frac{\rho}{p}, \frac{\gamma}{p}\Big\}$. This is the key remark for proving the next proposition.

\begin{prop}
\label{PropSummaryPropertiesGeneralCase}
\begin{enumerate}
   \item Suppose $V$ is Lipschitz-continuous and the $V_i$ are $(\gamma-1)$-Lipschitz. Then $\big(\mu_{ts}\big)_{0\leq s\leq t\leq T}$ is an approximate flow.  \vspace{0.1cm} 
   \item Suppose $V$ is $(1+\rho)$-Lipschitz, for some $\rho>\frac{p-[p]}{p}$, and the $V_i$ are $\gamma$-Lipschitz, then $\big(\mu_{ts}\big)_{0\leq s\leq t\leq T}$ is a $\mcC^1$-approximate flow. 
\end{enumerate}
\end{prop}

\medskip

\begin{Dem}
We first use formula \eqref{EqExactFormulaMuTs} to write
\begin{equation}
\mu_{tu}\big(\mu_{us}(x)\big) = \mu_{us}(x) + (t-u)V\big(\mu_{us}(x)\big) + \sum_I X^I_{tu} V_I\big(\mu_{us}(x)\big) +\ep^{\textrm{Id}\,;\,[p]}_{tu}\big(\mu_{us}(x)\big).
\end{equation}
We deal with the term $(t-u)V\big(\mu_{us}(x)\big)$ using \eqref{EqElementaryEstimateYu} and the Lipschitz character of $V$. The remainder $\ep^{\textrm{Id}\,;\,[p]}_{tu}\big(\mu_{us}(x)\big)$ has, under the assumptions of point (1), an infinite norm bounded above by $c\big(1+\|{\bfX}\|^\gamma\big)|t-u|^a$, and its derivative has infinite norm bounded above by $c\big(1+\|{\bfX}\|^\gamma\big)^2|t-u|^a$, under the assumptions of point (2), by the above remark and point (2) of proposition \ref{PropContinuityFlow} about the uniform Lipschitz size of the maps $\mu_{ba}$.

\ssk

To deal with the terms $X^I_{tu} V_I\big(\mu_{us}(x)\big)$, with $|I|=k$, we use formula \eqref{EqExactFormulaMuTs} with $n=\big([p]-k\big)$, to develop $V_I\big(\mu_{us}(x)\big)$. We have
\begin{equation*}
\begin{split}
&V_I\big(\mu_{us}(x)\big) = V_I(x) + (u-s)\big(VV_I\big)(x) + \sum_{j=1}^{[p]-k} \frac{1}{j!}\sum_{|I_1|+\dots+|I_j|\leq [p]} \left(\prod_{m=1}^j \Lambda^{I_m}_{us}\right) \big(V_{[I_j]}\cdots V_{[I_1]}V_I\big)(x)  + \ep^{V_I\,;\,{p}-k}_{us}(x) \\
&+\sum_{|I_1|+\dots+|I_{[p]-k}|\leq [p]} \left(\prod_{m=1}^{[p]-k} \Lambda^{I_m}_{us}\right) \int_0^1 \Big\{\Big(V_{[I_{[p]-k}]}\cdots V_{[I_1]}V_I\Big)\big(y_{s_{[p]-k}}\big)-\big(V_{[I_{[p]-k}]}\cdots V_{[I_1]}V_I\big)(x)\Big\}{\bf 1}_{\Delta_{[p]-k}}\,ds_{[p]-k}\dots ds_1.
\end{split}
\end{equation*}
We see using \eqref{EqElementaryEstimateYu} that $X^I_{tu}$ times the second line above has an infinite norm bounded above by $c\big(1+\|{\bfX}\|^\gamma\big)|t-s|^a$, and that its derivative has an infinite norm bounded above by $c\big(1+\|{\bfX}\|^\gamma\big)^2|t-s|^a$, under the regularity assumptions of point (2).

Writing 
$$
\mu_{us}(x) = x + (u-s)V(x) + \sum_{I} X^I_{us}V_I(x) + \ep^{\textrm{Id}\,;\,[p]}_{us}(x),
$$
it is then straightforward to use the identities $\exp\big({\bf \Lambda}_{us}\big)=\bfX_{us}$ and ${\bfX}_{ts} = {\bfX}_{us}{\bfX}_{tu}$, to see that
$$
\mu_{tu}\big(\mu_{us}(x)\big) = \mu_{ts}(x) + \ep_{ts}(x),
$$
with a remainder $\|\ep_{ts}\|_\infty \leq c\big(1+\|{\bfX}\|^\gamma\big)|t-s|^a$, under the regularity assumptions of point (1), with first derivative with infinite norm bounded above by $c\big(1+\|{\bfX}\|^\gamma\big)^2|t-s|^a$, under the assumptions of point (2).

It remains to check that the decomposition \eqref{EqH0} holds in that case, which is done by writing 
\begin{equation}
\label{EqDecoupageMu_ts}
\mu_{ts}(x) = x + \int_0^1\Big\{(t-s)V(y_r^x) + \sum_{|I|=[p]}\Lambda^I_{ts}V_{[I]}(y_r^x)\Big\}dr + \int_0^1 \sum_{|I|\leq [p]-1}\Lambda^I_{ts}V_{[I]}(y_r^x) dr
\end{equation}
(where we have emphasized the dependence of $y_r$ on its initial condition $x$ by an upper index $x$), and defining $A^{ts}_x$, as the differential of the function of $x$ defined by the first integral, and $B^{ts}_x$, by the differential of the function of $x$ defined by the second integral; with the notation of proposition \ref{PropUniformLipControl}, we have $\al= \frac{[p]}{p}$ here.
\end{Dem}

\medskip

The following is to be tought of as an analogue in the setting of flows of Davie's definition of a solution to a rough differential equation \cite{Davie}, as recalled in the introduction.

\medskip

\begin{defn}
\label{DefnGeneralRDESolution}
Let $2\leq p<\gamma< [p]+1$ be given. Let $V_1,\dots,V_\ell$ be $\mcC_b^{[p]}$-Lipschitz vector fields on $E$, and $\bfX$ be a H\"older weak geometric $p$-rough path. Write \emph{F} for $(V_1,\dots,V_\ell)$.  Let $V$ be a bounded Lipschitz continuous vector field on $E$. With the above notations,  a \textbf{\emph{flow}} $(\varphi_{ts}\,;\,0\leq s\leq t\leq T)$ is said to \textbf{\emph{solve the rough differential equation}}
\begin{equation}
\label{RDEGeneral}
d\varphi = Vdt + \textrm{\emph{F}}\,\bfX(dt)
\end{equation}
if there exists a constant $a>1$ independent of $\bfX$ and two possibly $\bfX$-dependent positive constants $\delta$ and $c$ such that
\begin{equation}
\label{DefnSolRDEGeneral}
\|\varphi_{ts}-\mu_{ts}\|_\infty \leq c\,|t-s|^a
\end{equation}
holds for all $0\leq s\leq t\leq T$ with $t-s\leq\delta$.
\end{defn}

The following well-posedness result follows directly from theorem \ref{CorMainWorkingResult} and proposition \ref{PropSummaryPropertiesGeneralCase}.

\begin{thm}
\label{ThmMainResultGeneral}
Suppose $V$ is $(1+\rho)$-Lipschitz, for some $\rho>\frac{p-[p]}{p}$, and the $V_i$ are $\gamma$-Lipschitz. Then the rough differential equation $d\varphi = Vdt + \textrm{\emph{F}}\,\bfX(dt)$ has a unique solution flow; it takes values in the space of uniformly Lipschitz continuous homeomorphisms of $E$ with uniformly Lipschitz continuous inverses, and depends continuously on $\bfX$.
\end{thm}

\ssk

\begin{Dem}
\noindent Note that any solution flow depends continuously on $(s,t)$ in the topology of uniform convergence on $E$ by proposition \ref{PropContinuityFlow} and \eqref{DefnSolRDEGeneral}.  Use the notations $c_1,c_2,c_3$ of section \ref{SectionConstructingFlows} for the constants appearing in the sewing lemma for flows, and write $a$ for $\min\Big\{1+\frac{\rho}{p},\frac{\gamma}{p}\Big\}$; it follows from the above estimates that we can choose 
$$
c_1=c_3=c\big(1+\|\bfX\|^\gamma\big),\quad c_2 = c\,c_3,
$$
so we have 
\begin{equation}
\label{EqApproxPhiMu}
\|\varphi_{ts}-\mu_{\pi_{ts}}\|_\infty \leq c\big(1+\|\bfX\|^\gamma\big)T\,\big|\pi_{ts}\big|^{a-1},
\end{equation}
for any partition $\pi_{ts}$ of $(s,t)\subset [0,T]$ of mesh $\big|\pi_{ts}\big|\leq \delta$, as a consequence of inequality \eqref{EqApproxVarphiMu}. As these bounds are uniform in $(s,t)$, and for $\bfX$ in a bounded set of the space of H\"older $p$-rough paths, and each $\mu_{\pi_{ts}}$ is a continuous function of $\bfX$, by proposition \ref{PropContinuityFlow}, the flow $\varphi$ depends continuously on $\big((s,t),\bfX\big)$.

To prove that $\varphi$ is a homeomorphism, note that, with the notations of section \ref{SectionConstructingFlows},
$$
\Big(\mu^{(n)}_{ts}\Big)^{-1} = \mu_{s_1s_0}^{-1}\circ\cdots\circ\mu_{s_{2^n}s_{2^n-1}}^{-1}
$$
can actually be written $\big(\mu^{(n)}_{ts}\big)^{-1} = \widetilde\mu_{s_{2^n}s_{2^n-1}}\circ\cdots\circ\widetilde\mu_{s_1s_0}$, for the time 1 map $\widetilde\mu$ associated with the rough path $\bfX_{t-\bullet}$. As $\widetilde\mu$ enjoys the same properties as $\mu$, the maps $\big(\mu^{(n)}_{ts}\big)^{-1}$ converge uniformly to some continuous map $\varphi_{ts}^{-1}$ which satisfies by construction $\varphi_{ts}\circ\varphi_{ts}^{-1} = \textrm{Id}$.

\ssk

As propositions \ref{PropUniformLipControl} and \ref{PropSummaryPropertiesGeneralCase} provide a uniform control of the Lipschitz norm of the maps $\mu_{ts}^{(n)}$, the limit maps $\varphi_{ts}$ also have Lipschitz norms controlled by the same quantity; the same holds for their inverses. We propagate this property from the set $\big\{(s,t)\in [0,T]^2\,;\,s\leq t, \; t-s\leq \delta\big\}$ to the whole of $\big\{(s,t)\in [0,T]^2\,;\,s\leq t\big\}$ using the flow property of $\varphi$.
\end{Dem}

\medskip

\begin{rems}
\label{Remarks}
\begin{enumerate}
   \item {\bf Friz-Victoir approach to rough paths.} The continuity of the solution flow \wrt the driving rough path ${\bf X}$ has the following consequence, which justifies the point of view adopted by Friz and Victoir in their works. Suppose the H{\"o}lder weak geometric $p$-rough path ${\bf X}$ is the limit in the rough path metric of the canonical H{\"o}lder weak geometric $p$-rough paths ${\bf X}^n$ associated with a smooth (or Lipschitz continuous) $E$-valued path $(x^n_t)_{0\leq t\leq T}$ through the data of its well-defined iterated integrals. The illustration to theorem \ref{CorMainWorkingResult} given at the end of section \ref{SectionConstructingFlows} shows that the solution flow $\varphi^n$ to the rough differential equation \eqref{RDEGeneral} with driving rough path ${\bf X}^n$ is the flow associated with the ordinary differential equation 
$$
\dot y_u = V(y_u)du+V_i(y_u)d(x^n_u)^i. 
$$
As $\|\varphi^n-\varphi\|_\infty=o_n(1)$, from the continuity of the solution flow \wrt the driving rough path, the flow $\varphi$ appears in that case as a uniform limit of the elementary flows $\varphi^n$. A H{\"o}lder weak geometric $p$-rough path with the above property is called a \textit{H{\"o}lder geometric $p$-rough path}; not all H{\"o}lder weak geometric $p$-rough path are H{\"o}lder geometric $p$-rough path \cite{FVGeometric}. 

Friz and Victoir obtain in \cite{FVEuler} convergence rates similar to \eqref{EqApproxPhiMu}, with a slightly better exponent, by a clever use of sub-Riemannian geometry in $G^{[p]}_\ell$, for rough differential equations without a drift. The above estimate for rough differential equations with a drift appears to be new.
\vspace{0.2cm}
	
	\item  \textbf{Time-inhomogeneous dynamics.} The above result have a straightforward generalization for a time dependent bounded drift $V(s;\cdot)$ which is Lipschitz continuous \wrt the time and $(1+\rho)$-Lipschitz \wrt the space variable, and time dependent bounded vector fields $V_i(s;\cdot)$ which are $\gamma$-Lipschitz \wrt the space variable and Lipschitz continuous \wrt time. We define in that case a $\mcC^1$-approximate flow by defining $\mu_{ts}$ as the time 1 map associated with the ordinary differential equation
$$
\dot y_u = (t-s)V(s;y_u) + \sum_{r=1}^{[p]}\sum_{I\in\llbracket 1,\ell\rrbracket^r} \Lambda_{ts}^{r,I} V_{[I]}(s;y_u),\quad 0\leq u\leq 1.	
$$
In particular, inequality \eqref{GralFundamentalEstimate} holds in that case, with $V(s;x)$ and $V_I(s;x)$ in place of $V(x)$ and $V_I(x)$. This framework will be useful in section \ref{SectionNonlinearRDE} for the study of some simple mean field stochastic rough differential equation. \vspace{0.2cm}
   
   \item {\bf Flows driven by (non-weak geometric) H{\"o}lder $p$-rough paths.} Note that the only place where we have used the fact that ${\bf X}$ takes values in the free nilpotent Lie group $G^{[p]}_\ell$ is in writing its logarithm as a sum of terms involving only the brackets of the basis vectors $e_1,\dots,e_\ell$. The above reasonnings are totally independent of that fact. Given any H\"older $p$-rough path $\bfX$, set $\log {\bf X}_{ts} = {\bf \Lambda}_{ts} = \sum_{r=1}^{[p]}\sum_{I\in\llbracket 1,\ell\rrbracket^r} \Lambda_{ts}^{r,I} e_I$, and define in that case $\mu_{ts}$ as the time one map associated with the ordinary differential equation 
   $$
   \dot y_u = (t-s)V(y_u) + \sum_{r=1}^{[p]}\sum_{I\in\llbracket 1,\ell\rrbracket^r} \Lambda_{ts}^{r,I}V_I(y_u).
   $$   
The above proofs show that theorem \ref{ThmMainResultGeneral} holds in that case as well, generalizing the results of Lyons \cite{Lyons97} and Lejay \cite{LejayGlobal}, who only deals with the case $2\leq p<3$, without a drift. The point of working with weak geometric rough paths is the fact that only brackets of vector fields appear in the definition of $\mu_{ts}$. These brackets are intrinsically defined vector fields which do not depend on any choice of a coordinate system, so are well-defined on a manifold. This is not the case of the vector fields $V_I$ which are needed in the above definition of $\mu_{ts}$ and are meaningless in a manifold setting. This fact justifies calling $G^{([p])}_\ell$-valued rough paths (weak) geometric rough paths. This way of solving rough differential equation driven by non-weak geometric H\"older $p$-rough paths offers an alternative to the use of branched rough paths \cite{GubinelliBranched}.
\end{enumerate}
\end{rems}

\ssk

\subsection{Taylor expansion of solution flows to rough differential equations}
\label{SectionTaylorExpansion}

It is a nice feature of our approach that Taylor expansion of solution flows to rough differential equations come almost for free, providing an alternative and simple proof of similar results due to Friz and Victoir \cite{FVEuler}. As in the classical case, the smoother the vector fields $V_i$ are, the better we can describe the solution flow to the driftless equation $d\varphi = \textrm{F}\bfX(dt)$. This takes the form of a refined version of the inequality $\big\|\varphi_{ts}-\mu_{ts}\big\|_\infty \leq c\,|t-s|^{\frac{\gamma}{p}}$. 

\medskip

Let $2\leq p$ and an integer $N\geq [p]+1$ be given. Lyons devised in his original theory of rough paths \cite{Lyons97} a fundamental mechanism which enables in particular to extend uniquely any H\"older weak geometric $p$-rough path to a very special $G_\ell^{(N)}$-valued map and provides the missing higher order iterated integrals needed to write a Taylor expansion of the solution flow to the equation $d\varphi = \textrm{F}\bfX(dt)$. We state it here under the form we need -- see \cite{FdlPM} for a simple proof of a refined version of Lyons' extension theorem.

\begin{thm}[Lyons' extension theorem]
Let $2\leq p$ be given and $\bfX$ be a H\"older weak geometric $p$-rough path. Let $N\geq [p]+1$ be an integer. Then there exists a unique ($G_\ell^{(N)}$-valued) H\"older weak geometric $N$-rough path path ${\bfY}$ extending $\bfX$ in the sense that $\pi_{[p] N}{\bfY}_{ts}={\bfX}_{ts}$, for all $0\leq s\leq t\leq T$. The map $\bfY$ is called the lift of $\bfX$ to $G_\ell^{(N)}$; it is a continuous function of $\bfX$.
\end{thm}

Given $0<\alpha<1$, set $\gamma=N+\alpha$, and suppose the vector fields $V_i$ are $\gamma$-Lipschitz. Given a H\"older weak geometric $p$-rough path $\bfX$, let $\bfY$ its lift to $G^{(N)}_\ell$, and ${\bf\Lambda} = 0\oplus\Lambda^1\oplus\Lambda^2\oplus\cdots\oplus \Lambda^N$ be the logarithm of its lift,  so we have $\exp {\bf \Lambda}_{ts} = {\bfY}_{ts}$, for all $0\leq s\leq t\leq T$. \vspace{0.1cm}

Denote by $\nu^{[N]}_{ts}$ the diffeomorphism of $E$ which associates to any $x\in E$ the value at time $1$ of the well-defined and unique solution of the ordinary differential equation
$$
\frac{d}{dr}y_u = \sum_{r=1}^N\sum_{I\in\llbracket 1,\ell\rrbracket^r} \Lambda^{r,I}_{ts}V_{[I]}(y_u), \quad 0\leq u\leq 1,
$$
with $y_0=x$. The proof of proposition \ref{PropEstimatesGeneralRDE} shows that the following awaited estimate holds.

\begin{prop}
\label{PropEstimteRhoN}
Let $2\leq p$ and $[p]+1\leq N<N+\alpha< N+1$ be given, with $N$ integer. Let $V_1,\dots,V_\ell$ be $(N+\alpha)$-Lipschitz vector fields on $E$. Let $\bfX=1\oplus X^1\oplus\cdots\oplus X^{[p]}$ be a H\"older weak geometric $p$-rough path with lift $\bfY$ to $G_\ell^{(N)}$. Then there exists a positive constant $c$, depending only on $M,\lambda,T$ and $\|\bfX\|$ and $f\in\mcC^\gamma$ \st 
\begin{equation}
\label{GeneralFundamentalEstimate}
\left\|f\circ\nu^{[N]}_{ts} - \left\{f + \sum_{r=1}^{N} \sum_{I\in\llbracket 1,\ell\rrbracket^r} Y^{r,I}_{ts} V_If\right\}\right\|_\infty \leq c\,|t-s|^{\frac{N+\alpha}{p}}
\end{equation}
holds for all $f\in\mcC^\gamma$. The maps $\nu^{[N]}_{ts}$ depend continuously on $\big((s,t),\bfX\big)$ in uniform topology.
\end{prop}

As in section \ref{SectionGeneralCase}, this fundamental estimate implies together with theorem \ref{CorMainWorkingResult} a well-posedness result.

\begin{thm}
\label{ThmVarphiN}
Under the hypotheses of proposition \ref{PropEstimteRhoN}, there exists a unique flow $\big(\varphi^{[N]}_{ts}\big)_{0\leq s\leq t\leq T}$ on $E$ for which there are two positive constants $\delta$ and $c$ \st 
\begin{equation}
\label{RDEInequality}
\Big\|\varphi^{[N]}_{ts}-\nu^{[N]}_{ts}\Big\|_\infty \leq c\,|t-s|^{\frac{N+\alpha}{p}}
\end{equation}
holds for all $0\leq s\leq t\leq T$, with $t-s\leq \delta$.
\end{thm}

Let now $\gamma\in \big(p,[p]+1\big)$. Since $\Big\|\nu^{[N]}_{ts}-\mu_{ts}\Big\|_\infty\leq c|t-s|^{\frac{\gamma}{p}}$, the identity
$$
\Big\|\varphi^{[N]}_{ts}-\mu_{ts}\Big\|_\infty\leq c|t-s|^{\frac{\gamma}{p}}
$$
holds for all $0\leq s\leq t\leq T$, sufficiently close, so $\varphi^{[N]}$ is the solution flow to the rough differential equation $d\varphi=\textrm{F}\bfX(dt)$. 

\begin{cor}[Euler estimates/Taylor expansion]
\label{CorEulerEstimates}
Let  $2\leq p$, an integer $[p]+1\leq N$ and $0<\alpha<1$ be given. Suppose the vector fields $V_1,\dots,V_\ell$ are $(N+\alpha)$-Lipschitz. Let $\bfX$ be a H\"older weak geometric $p$-rough path. Denote by $\varphi$ the unique solution flow to the driftless rough differential equation $d\varphi = \textrm{\emph{F}}\bfX(dt)$. Let $\bfY$ be the lift of $\bfX$ to $G^{(N)}_\ell$. Then there exists two positive constants $\delta$ and $c$ \st 
\begin{equation*}
\Big\|\varphi_{ts}-\nu^{[N]}_{ts}\Big\|_\infty \leq c\,|t-s|^{\frac{N+\alpha}{p}}
\end{equation*}
holds for all $0\leq s\leq t\leq T$ with $t-s\leq \delta$. We have in particular, 
   \begin{equation}
   \label{RoughTaylorEstimates}   
   \left\|\varphi_{ts} - \left\{\textrm{\emph{Id}} + \sum_{r=1}^N \sum_{I\in\llbracket 1,\ell\rrbracket^r} Y^{r,I}_{ts}V_I\right\}\right\|_\infty \leq c\,|t-s|^{\frac{N+\alpha}{p}}
   \end{equation}
   for all $0\leq s\leq t\leq T$, with $t-s\leq \delta$.

\end{cor}

\begin{rems}
\begin{enumerate}
   \item In a probabilistic context where $\bfX$ is the random realization of the Brownian rough path, the $Y^{k,I}_{ts}$ coincide almost-surely with the iterated Stratonovich integrals $\int_0^t \circ dB^{i_1}_{s_1}\otimes\cdots\otimes \circ dB^{i_k}_{s_k}$, for $I=(i_1,\dots,i_k)$, and \eqref{RoughTaylorEstimates} is a pathwise version of Azencott's celebrated stochastic Taylor formula -- see for instance \cite{Azencott}, \cite{BenArousTaylor} and \cite{Castell}. \vspace{0.1cm}
   \item Friz and Victoir \cite{FVEuler} proved similar estimates by a clever use of geodesic approximation in the free nilpotent Lie group $G^{(N)}_\ell$. \vspace{0.1cm}
   \item The above theorem provides a straightforward justification of Friz and Oberhauser's theorem about drift induced by perturbed driving signals; theorem 2 in \cite{FrizOberhauser}.
\end{enumerate}
\end{rems}

\section{Non-explosion under linear growth conditions on the vector fields}
\label{SectionLinearGrowth}

Let $1<\gamma$ be a non-integer real number. We say that a vector field $W$ is {\bf $\gamma$-Lipschitz with linear growth} if it is $\mcC^{[\gamma]}$, with bounded derivatives $W^{(1)},\cdots,W^{([\gamma])}$ and its $[\gamma]$'s derivative is a bounded $\big(\gamma-[\gamma]\big)$-H\"older map. We do not require that $W$ itself be bounded. We describe in that section how the arguments of sections \ref{SectionConstructingFlows} and \ref{SectionGeneralRDE} need to be amended to prove a well-posedness result for the rough differential equation on flows \eqref{BasicRDEquation}, under the relaxed assumption that the driving vector fields $V,V_1,\dots,V_\ell$ are sufficiently regular and have linear growth. We need for that purpose a suitable version of the sewing lemma for flows which applies in that setting. The main difficulty here is that the solution to the ordinary differential equation
$$
\dot y_u = (t-s)V(y_u) + \sum_{r=1}^{[p]}\sum_{I\in\llbracket 1,\ell\rrbracket^r} \Lambda_{ts}^{r,I}V_{[I]}(y_u)
$$
may explode in finite time since the vector fields $V_{[I]}$ have polynomial growth now. This tendancy is compensated by the small size of the coefficients $\Lambda_{ts}^{r,I}$, which ensures that the time $1$ map $\mu_{ts}$ will be well-defined on some ball of fixed radius provided $t-s$ is small enough. So we need to set up a framework where approximate flows are replaced by some kind of local approximate flows, not defined in the whole of $E$. This is what section \ref{SectionLocalApproxFlows} is about. 

\ssk

We use in that section the notation $B_R$ for the ball $\big\{x\in E\,;\,|x|\leq R\big\}$, and $c(R)$ for a constant depending only on $R$, whose precise value is unimportant.

\subsection{Local approximate flows with exponential growth}
\label{SectionLocalApproxFlows}

Let $\mu : \big(x,(s,t)\big)\mapsto \mu_{ts}(x)$ be a continuous map from 
$$
\bigcup_{x\in E} \{x\}\times\Big\{(s,t)\in [0,T]^2\,;\,0\leq t-s < ce^{-c_0|x|}\Big\}
$$ 
to $E$, where $c_0$ is a positive constant. For $0\leq s\leq t\leq T$ fixed, $\mu_{ts}(\cdot)$ is defined on the centered open ball with radius $\frac{-\ln(t-s)}{c_0}$, where we assume it is of class $\mcC^1$. The following assumption quantifies the fact that $\mu_{ts}$ is $\mcC^1$-close to the identity, for $(t-s)$ small. \vspace{0.1cm}

\begin{itemize}
   \item[{\bf H1'.}] {\it For all $R>0$ and $x\in B_R$, for all $0\leq s\leq t\leq T$ with $t-s\leq ce^{-c_oR}$, we have}
   \begin{equation*}
   \begin{split}
   &\textrm{{\bf a)}} \quad \mu_{ts}(x) = x + O_c\Big(|t-s|^{\frac{1}{2p}}\Big), \\
   &\textrm{{\bf b)}} \quad \big|D_x\mu_{ts}\big| \leq 1+ O_c\Big(|t-s|^{\frac{1}{2p}}\Big).   
  \end{split}
  \end{equation*}
\end{itemize}  
Fix $0<\ep\leq \frac{1}{2}$. It follows from {\bf H1' a)} that $\mu_{tu}\big(\mu_{us}(x)\big)$ makes sense for $x\in B_R$ and any $0\leq s\leq u\leq t\leq T$, with $t-s\leq ce^{-c_oR}$ and $\ep\leq\frac{u-s}{t-s}\leq 1-\ep$, provided $R$ is bigger than some radius $R_\ep$ depending only on $\ep$ and the data of the problem. The following assumption expresses the fact that $\mu$ is a "{\it local approximate flow}". \vspace{0.1cm}

\begin{itemize}
   \item[{\bf H2'.}] {\it There exists two positive constants $c_1$ and $a>1$ such that}
   \begin{equation}
   \label{EqMumuLinear}
   \big|\mu_{tu}\big(\mu_{us}(x)\big)-\mu_{ts}(x)\big| \leq ce^{c_1|x|}\,|t-s|^a
   \end{equation}   
   {\it holds for all $x\in B_R$ and all $0\leq s\leq u\leq t\leq T$ with $\ep\leq\frac{u-s}{t-s}\leq 1-\ep$, and $t-s\leq ce^{-c_0R}$.}
\end{itemize} 
Last, assume that 
\begin{itemize}
   \item[{\bf H3'.}] {\it there exists a function $\delta : (0,\infty)\rightarrow \RR_+$, no greater than $R\mapsto ce^{-c_oR}$, and for any $R>0$ an $R$-dependent contant $c_2(R)$, such that the restriction to the ball $B_R$ of the maps $\mu^n_{ts}(\cdot)$, for $n\geq 2$ and $t-s\leq \delta(R)$, are all Lipschitz continuous, with Lipschitz constants uniformly bounded above by $c_2(R)$.}  \vspace{0.1cm}
\end{itemize}

\begin{thm}
\label{ThmSewingLemmaFlowsLinear}
A map $\mu$ satisfying assumptions \emph{{\bf H1',H2',H3'}}, defines a unique flow $\big(\varphi_{ts}\big)_{0\leq s\leq t\leq T}$ on $E$ to which one can associate a function $\varep : (0,\infty)\rightarrow \RR_+$, such that 
   \begin{equation}
   \label{EqApproxVarphiLinear}
   \big|\varphi_{ts}(x)-\mu_{ts}(x)\big| \leq c(R)\,|t-s|^a
   \end{equation}
holds for all $x\in B_R$, and all $0\leq s\leq u\leq t\leq T$, with $t-s\leq \varep(R)$, for all $R>0$.
\end{thm}

A map enjoying property \eqref{EqMumuLinear} is called a {\bf local approximate flow with exponential growth}. The proof of theorem \ref{ThmSewingLemmaFlowsLinear} requires the following version of lemma \ref{LemImprovedEqMumu}. Set $c_\ep = \underset{\ep\leq \beta\leq 1-\ep}{\max}\;\big\{\beta^a+(1-\beta)^a\big\}<1$, and 
$$
L = \frac{2c_1}{1-c_\ep}\vee\Big(a-\frac{1}{2p}\Big)c_0.
$$
Set $\delta_1(R) = c\exp\Big(\frac{-1}{a}\big(L+\frac{c_0}{2p}\big)R\Big)\leq ce^{-c_0R}$. This function satisfies for any $0\leq s\leq u\leq t\leq T$ with $t-s\leq \delta_1(R)$, and $\ep \leq \frac{u-s}{t-s}\leq 1-\ep$, the inequality
\begin{equation}
\label{EqDelta1}
\delta_1\Big(R+ce^{-\frac{c_0}{2p}R}\Big) \geq (1-\ep)\delta_1(R) \geq t-u.
\end{equation}

\begin{lem}
\label{LemEqMumuLinear}
Let $\mu$ be a local approximate flow with exponential growth satisfying assumptions \emph{{\bf H1'}}. Given $\ep>0$, there exists a constant $R_\ep$ depending on the data of the problem such that for any $R\geq R_\ep$, any $0\leq s\leq t\leq T$, with $t-s\leq \delta_1(R)$, and any special partition $\pi_{ts}$ of $(s,t)$ of type $\ep$, the map $\mu_{\pi_{ts}}(\cdot)$ is well-defined on $B_R$, and we have for any $x\in B_R$ 
\begin{equation}
\label{EqImprovedEqMumuLinear}
\big|\mu_{\pi_{ts}}(x)-\mu_{ts}(x)\big| \leq Le^{L|x|}\,|t-s|^a.
\end{equation}
\end{lem}

\ssk

\begin{Dem}
Take $0\leq s\leq t\leq T$, with $t-s\leq \delta_1(R)$. We proceed by induction on the number of sub-intervals of the partition $\pi_{ts}$, as in the proof of lemma \ref{LemImprovedEqMumu}. With the notations adopted there, the induction hypothesis and assumption {\bf H1' a)} tell us that $\mu_{\pi_{us}}(x)$ is well-defined and has norm bounded above by 
$$
\Big|\mu_{\pi_{us}}(x)\Big| \leq R + c|u-s|^{\frac{1}{2p}} + Le^{LR}|u-s|^a \leq R + ce^{-\frac{c_0}{2p}R},
$$
so \eqref{EqDelta1} ensures that $\mu_{tu}\big(\mu_{us}(x)\big)$ is well-defined for $R$ big enough. For such radii, we have from the decomposition given in \eqref{EqFirstEstimate}
\begin{equation*}
\big|\mu_{\pi_{ts}}(x) - \mu_{ts}(x)\big| \leq Le^{L|\mu_{\pi_{us}(x)}|}|t-u|^a + \Big(1+ce^{-\frac{c_0R}{2p}}\Big)Le^{L|x|}|u-s|^a + c_1e^{c_1|x|}|t-s|^a,
\end{equation*}
by using the induction hypothesis and assumptions {\bf H1' a,b)} to get the first two terms on the right hand side, and \eqref{EqMumuLinear} to get the third term. 
So, setting $o_R(1) = ce^{-\frac{c_0}{2p}R}$ and $u-s =\beta (t-s)$, with $\ep\leq \beta\leq 1-\ep$, we have the inequality
\begin{equation*}
\begin{split}
\big|\mu_{\pi_{ts}}(x) - \mu_{ts}(x)\big| &\leq \Big\{\Big(e^{o_R(1)}(1-\beta)^a +\Big(1+o_R(1)\big)\beta^a\Big)L+c_1\Big\} e^{L|x|} |t-s|^a \\
&\leq \big(e^{o_R(1)}c_\ep L+c_1\big)e^{L|x|} |t-s|^a.
\end{split}
\end{equation*}
One has $e^{o_R(1)}c_\ep L+c_1 \leq L$, provided $R$ is bigger than some explicit radius $R_\ep$ depending only on $\ep$.
\end{Dem}

\bigskip

\begin{DemSewingLemmaLinear}
We follow the proof of theorem \ref{ThmSewingLemmaFlows}, indicating the points that need to be amended.

\ssk

\noindent {\bf a) Existence.} Recall the definition of $\mu^{(n)}_{ts}$ given in \eqref{EqMu(n)}, for dyadic times $0\leq s\leq t\leq T$. Given $0\leq i\leq N(n)-1$ and $1\leq j\leq N(n)$, set 
\begin{equation*}
\begin{split}
&f_i = \mu_{s_{i+1}s_i+2^{-n-1}}\circ\mu_{s_i+2^{-n-1}s_i},\quad\textrm{and }\quad f_{[i]} = f_i\circ\cdots\circ f_0, \\
&g_j = \mu_{s_{j+1}s_j}\quad\textrm{and }\quad g^{[j]} = g_{N(n)-1}\circ\cdots\circ g_{N(n)-j}.
\end{split}
\end{equation*}
Set $\varep(R) = \delta(R)\wedge\delta_1(R)$, for all $R>0$, and take $s$ and $t$ such that $0\leq t-s\leq \varep(R)$. Then all the maps $f_{[i]}$ satisfy for all $x\in B_R$ the inequality
\begin{equation}
\label{EqEstimateF[i]}
\big|f_{[i]}(x)\big| \leq R+1.
\end{equation}
Also, as the maps $g^{[j]}$ are of the form $\mu^m_{vu}$, for some $m\geq 1$ and $s\leq u\leq v\leq t$, there exists for any positive $R'$ a constant $c_2(R')$ such that the restrictions to $B_{R'}$ of all the maps $g^{[j]}$ have Lipschitz constants uniformly bounded above by $c_2(R')$, by assumption {\bf H3'}. On the other hand, given $R>0$, it is a consequence of lemma \ref{LemEqMumuLinear} and assumption {\bf H1'} that all the maps $f_{N(n)-1-i}\circ f_{[N(n)-2-i]}$ and $g_{N(n)-1-i}\circ f_{[N(n)-2-i]}$ send the ball $B_R$ into a fixed ball $B_{R'}$. So, for $x\in B_R$, we have
\begin{equation}
\begin{split}
\Big|\mu^{(n+1)}_{ts}(x)-\mu^{(n)}_{ts}(x)\Big| &\leq c_2(R')\sum_{i=0}^{N(n)-1} \left|\big(f_{N(n)-1-i}-g_{N(n)-1-i}\big)\big(f_{[N(n)-2-i]}(x)\big)\right| \\
&\leq c_1c_2(R') \sum_{i=0}^{N(n)-1} e^{c_1\big|f_{[N(n)-2-i]}(x)\big|}\,\big|s_{N(n)-i}-s_{N(n)-i-1}\big|^a \\
&\leq c_1c_2(R')e^{(R+1)c_1} \sum_{i=0}^{N(n)-1} \big|s_{i+1}-s_i\big|^a;
\end{split}
\end{equation}
the second inequality comes from assumption {\bf H2'} and the third from \eqref{EqEstimateF[i]}. The existence of the flow $\varphi$ follows from that identity as in the proof of theorem \ref{ThmSewingLemmaFlows}.

\ssk

\noindent {\bf b) Uniqueness.} One proceeds as in the proof of theorem \ref{ThmSewingLemmaFlows} by replacing the uniform estimate $\psi_{ts} = \mu_{ts} + O_c\big(|t-s|^a\big)$, by the local estimate $\psi_{vu}(x) = \mu_{vu}(x) + O_c\big(|v-u|^a\big)$, which holds for all  $s\leq u\leq v\leq t$ sufficiently close, all $|x|\leq R$ and an $R$-dependent constant $c$.
\end{DemSewingLemmaLinear}

\medskip

\noindent Using lemma \ref{LemEqMumuLinear} instead of lemma \ref{LemImprovedEqMumu}, a direct adaptation of the proof of proposition \ref{PropUniformLipControl} provides its local counterpart, which we can use as a sufficient condition for condition {\bf H3'} to hold.

\begin{prop}[Uniform Lipschitz controls of local approximate flows with exponential growth]
\label{PropUniformControlsLinearGrowth}
Let $\alpha$ and $\rho$ be positive constants, with $0<a(1-\rho)<\alpha<1$, and 
\begin{itemize}
   \item[{\bf H0'.}] let $\mu$ be a local approximate flow with exponential growth on $E$, made up of $(1+\rho)$-Lipschitz locally defined maps, such that one can write
   $$
   \quad D_x\mu_{ts} = \textrm{\emph{Id}} + A^{ts}_x + B^{ts}_x
   $$ 
for some $\textrm{\emph{L}}(E)$-valued $\rho$-Lipschitz maps $A^{ts}$ whose restrictions to every ball $B(R)$ has $\rho$-Lipschitz norm bounded above by $c(R)|t-s|^\alpha$, and some maps $\textrm{\emph{L}}(E)$-valued maps $B^{ts}$, of class $\mcC^1$, whose restrictions to any ball $B(R)$ have $\mcC^1$-norm bounded above by $o_{c(R)|t-s|^{1/p}}(1)$. 
\end{itemize}
If there exists a positive constant $c_3$ such that we have for all $R>0$ and all $x\in B_R$
\begin{equation}
\label{EqUniformLipControlsLinear}
\big|D_x(\mu_{tu}\circ\mu_{us}) - D_x\mu_{ts}\big| \leq ce^{c_3|x|}\,|t-s|^a
\end{equation}
for all $0\leq s\leq u\leq t\leq T$ with $\frac{1}{3} \leq \frac{u-s}{t-s}\leq \frac{2}{3}$ and $t-s \leq \delta_1(R)$, then for every constant $c_4>c_3$, there exists some positive constant $R'_\ep$, and  a function $\delta : \big[R'_\ep,\infty\big) \rightarrow \RR_+$, such that 
$$
\big|D_x\mu^n_{ts} - D_x\mu_{ts}\big| \leq ce^{c_4R}\,|t-s|^a
$$
holds for all $n\geq 2$ and all $x\in B_R$, whenever $t-s\leq \delta(R)$, and $R\geq R'_\ep$.
\end{prop}

Note that the regularity condition {\bf H1' b)} holds under condition {\bf H0'}. A local approximate flow with exponential growth satisfying the regularity conditions {\bf H0'} and {\bf H1' a)}, and condition \eqref{EqUniformLipControlsLinear} is said to be a {\bf local $\mcC^1$-approximate flow with exponential growth}. Theorem \ref{ThmSewingLemmaFlowsLinear} and proposition \ref{PropUniformControlsLinearGrowth} provide a machinery for constructing flows from local $\mcC^1$-approximate flows with exponential growth. 

\begin{thm}
\label{CorMainWorkingResultLinear}
A local $\mcC^1$-approximate flow with exponential growth defines a unique flow $(\varphi_{ts})_{0\leq s\leq t\leq T}$ on $E$ to which one can associate for all positive $R$ some $R$-dependent constants $c$ and $\delta$ such that
$$
\big|\varphi_{ts}(x)-\mu_{ts}(x)\big| \leq c|t-s|^a
$$
holds for all $x\in B_R$, whenever $t-s\leq \delta$.
\end{thm}

\subsection{Flows driven by weak geometric rough paths and vector fields with linear growth}

Let us go back to the setting and notations of section \ref{SectionGeneralCase} and consider the ordinary differential equation 
\begin{equation}
\label{EqODEApproxRDE}
\dot y_u = (t-s)V(y_u) + \sum_{r=1}^{[p]}\sum_{I\in\llbracket 1,\ell\rrbracket^r} \Lambda_{ts}^{r,I}V_{[I]}(y_u).
\end{equation}
Assume the vector field $V$ is $(1+\rho)$-Lipschitz with linear growth, with $\rho>\frac{p-[p]}{p}$, and the $V_i$ are $\gamma$-Lipschitz with linear growth. Given $x\in E$, the maximal solution to equation \eqref{EqODEApproxRDE} may well explode in finite time, due to the polynomial growth of the vector fields $V_{[I]}$. The function $\big|y_u\big|$ satisfies however an integral equation of the form
$$
\big|y_u\big| \leq |x| + a_{ts} + b_{ts}\int_0^u\big|y_r\big|^{[p]}dr,
$$
with positive coefficients and $b_{ts}\leq c|t-s|^{1/p}$, which guarantees that $y_u$ cannot explode before time $\frac{c}{|x^{[p]-1}||t-s|^{1/p}}$. So the time $1$ map $\mu_{ts}$ is well-defined on the set 
$$
\bigcup_{x\in E} \{x\}\times\Big\{(s,t)\in [0,T]^2\,;\,0\leq t-s < ce^{-c_0|x|}\Big\}, 
$$
for a good choice of constants $c, c_0$. For $x\in B_R$ and $0\leq t-s<ce^{-c_0R}$, one gets from the integral version of equation \eqref{EqODEApproxRDE} the elementary estimate
\begin{equation}
\label{EqEstimatey_u}
\big|y_u - x\big| \leq c\Big(1+\|{\bfX}\|^{[p]}\Big) \big(1+|x|\big) \exp\Big(c\Big(1+\|{\bfX}\|^{[p]}\Big)|t-s|^{\frac{1}{p}}\Big)|t-s|^\frac{1}{p}.
\end{equation}
As $R|t-s|^{\frac{1}{2p}}\leq cRe^{-\frac{c_0R}{2p}}\leq c$, condition {\bf H1' a)} holds as a consequence of the above estimate. Condition {\bf H1' b)} holds for the same reason, and it follows from \eqref{EqEstimatey_u} and formula \eqref{EqExactFormulaMuTs} that one has
$$
\Big|\ep^{f\,;\,[p]}_{ts}(x)\Big| \leq c\big(1+\|{\bfX}\|^\gamma\big)\|f\|_\gamma e^{c\|{\bfX}\|^{[p]}|x|} |t-s|^{\frac{\gamma}{p}}.
$$
Similar estimates hold for the remainders $\ep^{f\,;n}_{ts}$, for $1\leq n\leq [p]-1$, and their first derivatives, with some exponent $a>1$ in place of $\gamma/p$, as in the proof of proposition \ref{PropSummaryPropertiesGeneralCase}. So it is straightforward to follow the proof of proposition \ref{PropSummaryPropertiesGeneralCase} and see that $\mu$ is a local $\mcC^1$-approximate flow for which one can choose for $ce^{c_1|x|}$ and $ce^{c_3|x|}$ in \eqref{EqMumuLinear} and \eqref{EqUniformLipControlsLinear}, respectively, an expression of the form $c\big(1+\|{\bfX}\|^\gamma\big)e^{c\|{\bfX}\|^{[p]}|x|}$. The decomposition of $\mu_{ts}$ given in \eqref{EqDecoupageMu_ts} provides the definition of $A^{ts}$ and $B^{ts}$ used in checking that assumption {\bf H0'} holds.

\ssk

The flow $\big(\varphi_{ts}\big)_{0\leq s\leq t\leq T}$ uniquely determined by $\mu$, by theorem \ref{CorMainWorkingResultLinear}, is called the solution flow to the rough differential equation 
$$
d\varphi = Vdt+\textrm{F}{\bfX}(dt),
$$ 
on flows. It follows in particular from this well-posedness result that there is no explosion, which generalizes to our infinite dimensional setting the corresponding results of \cite{LejayGlobal} and \cite{LejayNonExplosion}, which only treat the case $2\leq p<3$ without a drift, and the finite dimensional analogue result of Friz-Victoir \cite{FVBook}, who work with geometric rough paths.

\section{Paths driven by rough paths}
\label{SectionPathRDE}

The results of section \ref{SectionGeneralRDE} provide an easy proof of sharp existence and well-posedness results for point dynamics driven by a H{\"o}lder weak geometric rough path. We adopt the notations of the preceeding sections.

\begin{defn}
Let $2< p<\gamma\leq [p]+1$ be given. Let $V_1,\dots,V_\ell$ be vector fields on $E$ of class $\mcC^{[p]}$, and $\bfX$ be a H\"older weak geometric $p$-rough path on the time interval $[0,T]$. Let $V$ be a Lipschitz continuous vector field. A \textbf{path} $(z_s\,;\,0\leq s\leq T)$ is said to \emph{\textbf{solve the rough differential equation}}
\begin{equation}
\label{PathRDEGeneral}
dz = Vdt + \textrm{\emph{F}}\,\bfX(dt)
\end{equation}
if there exists a constant $a>1$ independent of $\bfX$, and two possibly $\bfX$-dependent positive constants $\delta$ and $c$, such that
\begin{equation}
\label{DefnSolRDEGeneral}
\big|z_t-\mu_{ts}(z_s)\big| \leq c\,|t-s|^a
\end{equation}
holds for all $0\leq s\leq t\leq T$, with $t-s\leq\delta$.
\end{defn}

\ssk

Recall the compactness assumption {\bf (C)} on the vector fields $V_i$ stated in section \ref{SectionGeneralCase} and the fact that if the $V_i$ are only $(\gamma-1)$-Lipschitz, we denote by $\mu_{ts}(x)$ the time 1 value of any solution to equation \eqref{DeterministicODEGeneralCase} started from $x$, whose existence is garanteed by assumption {\bf (C)}.

\medskip

\begin{thm}[Existence/well-posedness]
\label{ThmExistenceWellPosednessRDEPath}
\hfill
\begin{itemize}
   \item[{\bf a)}] Let $V$ be a bounded Lipschitz continuous vector field that maps balls into compact sets, and suppose that the $(\gamma-1)$-Lipschitz vector fields $V_i$ satisfy the compactness assumption \emph{{\bf (C)}}. Then the rough differential equation \eqref{PathRDEGeneral} has a solution path.  \vspace{0.1cm}
   \item[{\bf b)}] This solution path is unique if the $V_i$ are $\gamma$-Lipschitz and $V$ is $(1+\rho)$-Lipschitz, for some $\rho>\frac{p-[p]}{p}$; it is in that case a continuous function of $\bfX$.	
\end{itemize}
\end{thm}

\ssk

\noindent Note that we do not make any compactness assumption on $V$ or the $V_i$ in part b) of the above statement. This part of theorem \ref{ThmExistenceWellPosednessRDEPath} is usually called {\bf Lyons' universal limit theorem} \cite{Lyons97}. 

\medskip

\begin{Dem}
\noindent \textbf{a) Existence.} Given $\ep$ and $t\in [0,T]$, with $k\ep\leq t<(k+1)\ep$, set 
$$
z_t^\ep = \big(\mu_{t\,k\ep}\circ\mu_{k\ep\,(k-1)\ep}\circ\cdots\circ\mu_{\ep\,0}\big)(z_0).
$$
We prove that the paths $(z^\ep_t)_{0\leq t\leq T}$ have a uniformly converging subsequence by showing that they form an equicontinuous family and that the set $\big\{\mu^\ep_t\big\}_{0<\ep\leq T}$ is precompact for any fixed $t\in [0,T]$.

\ssk

First, for $k\ep\leq s<t\leq (k+1)\ep$, it is elementary to see that $z^\ep_t-z^\ep_s=O_c\Big\{|t-s|^{\frac{1}{p}}\Big\}$, for some constant $c$ independent of $\ep$. Now, given $0\leq s<t\leq T$, with $k\ep\leq s<(k+1)\ep$, and $\ell\ep\leq t<(\ell+1)\ep$, with $k+1\leq\ell$, one has
$$
z^\ep_t = \mu_{t\,\ell\ep}\circ\left(\overset{\ell-1}{\underset{n=k+1}{\bigcirc}}\mu_{(n+1)\ep\,n\ep}\right)\circ\mu_{(k+1)\ep\,s}\big(z_s^\ep\big).
$$
As every partition of any interval into sub-intervals of equal length is of special type $\frac{1}{2}$, one has by lemma \ref{LemImprovedEqMumu} and proposition \ref{PropContinuityFlow}
$$
\overset{\ell-1}{\underset{n=k+1}{\bigcirc}}\mu_{(n+1)\ep\,n\ep} = \mu_{\ell\ep\,(k+1)\ep} + O_L\Big(\big|(k-\ell)\ep\big|^\frac{\gamma}{p}\Big)
$$
provided $t-s\leq \delta$. Since the  maps $\mu_{ab}$ are uniformly Lipschitz continuous, the approximate flow property of $\mu$ gives
\begin{equation}
\label{AlmostRDESolution}
z^\ep_t = \mu_{ts}(z^\ep_s) + O_{L+c}\Big(|t-s|^\frac{\gamma}{p}\Big).
\end{equation}
The equicontinuity of the family of paths $(z^\ep_\bullet)_{0<\ep\leq 1}$ follows from the above inequality and the fact that $\big\|\mu_{ts}-\textrm{Id}\big\|_\infty\leq c|t-s|^{\frac{1}{p}}$. 

\ssk

Fix $t\in [0,T]$. To see that the set $\big\{\mu^\ep_t\big\}_{0<\ep\leq T}$ is pre-compact, note first that for any open ball $B$ of $E$, the set $\big\{\mu_{ba}(z)-z\,;\,z\in B,\,0\leq a\leq b\leq T\big\}$ is compact. Indeed, as 
$$
\mu_{ba}(z)-z = (b-a)\int_0^1 V(y_u)du + \Lambda^I_{ba}\int_0^1V_{[I]}(y_u)du,
$$
and $(y_u)_{0\leq u\leq 1}$ remains in a fixed ball for $y_0$ ranging in $B$, the above integrals belong to a fixed compact convex set independent of $a$ and $b$, by the compactness assumption on $V$ and assumption {\bf (C)} on the $V_i$, and a well-known theorem of Mazur on convex hulls of compact subsets of Banach spaces. It follows from this fact and \eqref{AlmostRDESolution} that we have for any $N\geq 1$ and $s_i=i\frac{t}{N}$, 
\begin{equation*}
\begin{split}
z^\ep_t - z_0 &= \sum_{i=0}^{N-1}\Big\{\big\{\mu_{s_{i+1}s_i}\big(z_{s_i}^\ep\big)-z_{s_i}^\ep\big\} + O_{L+c}\Big((t/N)^{\frac{\gamma}
{p}}\Big)\Big\} \\
             &= o_N(1) + \sum_{i=0}^{N-1} \big\{\mu_{s_{i+1}s_i}\big(z_{s_i}^\ep\big)-z_{s_i}^\ep\big\}.
\end{split}
\end{equation*}
So there exists a compact set $K(N)$ depending only on $N$ such that 
$$
\big\{\mu^\ep_t\big\}_{0<\ep\leq T} \subset\big\{K(N)+o_N(1)\big\},
$$
which implies the pre-compactness of the set $\big\{\mu^\ep_t\big\}_{0<\ep\leq T}$.

\ssk

Ascoli's theorem ensures as a consequence the pre-compactness of the set of paths $\big\{\mu^\ep_\bullet\big\}_{0<\ep\leq T}$ in the uniform topology. We see that any limit path satisfies \eqref{DefnSolRDEGeneral} for all $0\leq s\leq t\leq T$ with $t-s\leq\delta$, by sending $\ep$ to $0$.

\medskip

\noindent \textbf{b) Uniqueness.} Set $a=\min\Big(1+\frac{\rho}{p},\frac{\gamma}{p}\Big)$. To prove uniqueness of the solution when the $V_i$ are $\gamma$-Lipschitz and $V$ is $(1+\rho)$-Lipschitz, note first that $z_\bullet = \varphi_{\bullet 0}(z_0)$ is a solution, and remark that since it satisfies \eqref{DefnSolRDEGeneral} we have for any $k\geq 1$
\begin{equation*}
\begin{split}
z_{k\ep}-z^\ep_{k\ep} &= \mu_{k\ep\,(k-1)\ep}(z_{(k-1)\ep}) + O_c\big(\ep^a\big) - \mu_{k\ep\,(k-1)\ep}\big(z^\ep_{(k-1)\ep}\big) \\
                                   &= \mu_{k\ep\,(k-1)\ep}\Big( \mu_{(k-1)\ep\,(k-2)\ep}(z_{(k-2)\ep}) + O_c\big(\ep^a\big)\Big) + O_c\big(\ep^a\big) - \mu_{k\ep\,(k-1)\ep}\big(z^\ep_{(k-1)\ep}\big) \\
                                   &= \big(\mu_{k\ep\,(k-1)\ep}\circ \mu_{(k-1)\ep\,(k-2)\ep}\big)(z_{(k-2)\ep}) + O_c\big(\ep^a\big) + O_c\big(\ep^a\big) - \mu_{k\ep\,(k-1)\ep}\big(z^\ep_{(k-1)\ep}\big) 
\end{split}
\end{equation*}
as the $\mu_{k\ep\,(k-1)\ep}$ are Lipschitz continuous, uniformly \wrt $k$ and $\epsilon\leq 1$. As all the $\mu_{k\ep\,(k-1)}\circ\cdots\circ\mu_{\ell\ep\,(\ell-1)\ep}$ are Lipschitz continuous, with a \emph{uniform} upper bound for their Lipschitz constants for any $k\leq\ell$ and $\epsilon\leq 1$, by proposition \ref{PropSummaryPropertiesGeneralCase} (3), for which we need the $V_i$ to be $\gamma$-Lipschitz, and proposition \ref{PropUniformLipControl}, a basic induction gives
\begin{equation*}
\begin{split}
z_{k\ep}-z^\ep_{k\ep} &= \big(\mu_{k\ep\,(k-1)}\circ\cdots\circ\mu_{\ep\,0}\big)(z_0) + O_c\big(\ep^{-1}\big)\,O_c\big(\ep^a\big) -  \big(\mu_{k\ep\,(k-1)}\circ\cdots\circ\mu_{\ep\,0}\big)(z_0) \\
&= O_c\big(\ep^{-1}\big)\,O_c\big(\ep^a\big) = o_\ep(1).
\end{split}
\end{equation*}
We conclude by sending $\ep$ to $0$ along a subsequence for which $z^\ep_\bullet$ converges uniformly.
\end{Dem}

\begin{rems}
\begin{enumerate}
   \item 
   One can reformulate condition \eqref{DefnSolRDEGeneral} under the form
   $$
   \left|z_t-\left\{z_s+(t-s)V(z_s)+\sum_{r=1}^{[p]}\sum_{I\in\llbracket 1,\ell\rrbracket^r} X^{r,I}_{ts}V_I(z_s)\right\}\right| \leq c|t-s|^a 
   $$
   using proposition \ref{PropEstimatesGeneralRDE}. This form of the definition of a solution to a rough differential equation is due to Davie \cite{Davie} in the case $2\leq p<3$; it was extended to the general setting of H\"older geometric $p$-rough paths by Friz-Victoir \cite{FVEuler} under a different form. \vspace{0.1cm}
   \item The existence statement a) was first proved in a finite dimensional setting by Davie \cite{Davie}. Caruana \cite{Caruana} obtained an existence result in a Banach space setting, working with a full rough differential equation under a similar compactness assumption as {\bf (C)}, but for a rough differential equation with no drift term. As full rough differential equations enter the framework of rough differential equations driven by vector fields with linear growth, one recovers and extends in a simpler and (much) shorter way Caruana's result using lemma \ref{LemEqMumuLinear} instead of lemma \ref{LemImprovedEqMumu} in the above existence proof. 
\end{enumerate} 
\end{rems}

\section{An illustration: Mean field stochastic rough differential equations}
\label{SectionNonlinearRDE}

We show in this section how the results of sections \ref{SectionGeneralRDE} can be used to study some simple mean field stochastic rough differential equations. This kind of dynamics pops in naturally in the study of the large population limit of some classes of interacting random evolutions. The interaction holds through the dependence of the local characteristics of the random motion of each particle on the empirical measure of the whole family of particles. In a diffusion setting, each particle $i$ would satisfy a stochastic differential equation of the form
$$
dx^{(i)}_t = b\Big(x^{(i)}_t,\mu_t^{N}\Big)dt + \sigma\Big(x^{(i)}_t,\mu_t^{N}\Big)dB^{(i)}_t,
$$
where $\mu_t^N = \frac{1}{N}\sum_{k=1}^N\delta_{x^{(k)}_t}$. A large industry has been devoted to showing that the limit distribution in paths space of a typical particle of the system when $N$ tends to infinity has a dynamics of the form
\begin{equation}
\label{MeanFieldSDE}
dx_t = b\big(x_t,\mcL(x_t)\big)dt+\sigma\big(x_t,\mcL(x_t)\big)dB_t,
\end{equation}
where $\mcL(x_t)$ stands for the law of $x_t$. Theorem \ref{ThmMeanFieldRDE} provides a well-posedness result for such a limit equation, in the context of rough differential equations. As emphasized in \cite{CassLyons}, almost all the works in this area are set in the framework of a filtered probability space and rely crucially on some martingale arguments. On the other hand, the increasing importance of non-semi-martingale processes, like fractional Brownian motion, makes it desirable to have some more flexible tools to investigate equation \eqref{MeanFieldSDE} in such contexts. The theory of rough paths developped above provides a nice framework for that.

\medskip

A few notations are needed to set the probem. Given $2\leq p<\gamma\leq [p]+1$, we equip the set $\mcM_1\big(\RR^d\big)$ of probability measures on $\RR^d$ with the metric induced by its embedding in the dual of $\mcC^\gamma\big(\RR^d\big)$:
$$
\textrm{d}(\mu,\nu) = \sup\,\Big\{(g,\mu)-(g,\nu)\,;\,g\in\mcC^\gamma\big(\RR^d\big),\,\|g\|_\gamma \leq 1\Big\}.
$$
This metric topology is stronger than the weak convergence topology. Given any positive constant $m$, note that the set $\textrm{Lip}(m)$ of Lipschitz continuous paths from $[0,T]$ to $\big(\mcM_1\big(\RR^d\big),\textrm{d}\big)$, with Lipschitz constant no greater than $m$, is closed under the norm of uniform convergence, so it is a Banach space. 

\ssk

Fix $T>0$ and suppose $\bfX$ is a random variable defined on some probability space $(\Omega,\mcF,\PP)$, with values in the set of H\"older weak geometric $p$-rough paths over $\RR^\ell$, on the time interval $[0,T]$; write ${\bfX} = 1\oplus X^1\oplus\cdots\oplus X^{[p]}$, and $(\mcF_t)_{0\leq t\leq T}$ for the filtration generated by $\bfX$. As above, F stands for a collection $(V_1,\dots,V_\ell)$ of regular enough vector fields on $\RR^d$. Given a Lipschitz continuous path $\mcP = (P_t)_{0\leq t\leq T}$ in $\big(\mcM_1\big(\RR^d\big),\textrm{d}\big)$ and $\omega\in\Omega$, denote by $x_\bullet(\omega)$ the unique solution to the rough differential equation on paths
\begin{equation}
\label{NonlinearRDE}
dx_t = V(x_t,P_t)dt + \textrm{F}(x_t,P_t)\,\bfX(dt),
\end{equation}
where $x_0$ may be an integrable random variable independent of ${\bf X}$. Denote by $\Phi(\mcP)_t$ the law of $x_t$; a solution $x_\bullet$ to \eqref{NonlinearRDE} for which $\Phi(\mcP)_t = P_t$ is said to be a solution of the {\bf nonlinear rough differential equation} (or mean field stochastic rough differential equation)
\begin{equation*}
dx_t = V\big(x_t,\mcL(x_t)\big)dt + \textrm{F}\big(x_t,\mcL(x_t)\big)\,\bfX(dt).
\end{equation*}
Theorem \ref{ThmMeanFieldRDE} below provides conditions on $V$,F and ${\bf X}$ under which existence and uniqueness for solutions of this equation can be proved. Note that the main assumption \eqref{ConditionLawX} below holds for instance for the Brownian rough path, the rough path above an Orsntein-Uhlenbeck process, or fractional Brownian motion with Hurst index no smaller than $\frac{1}{2}$. 

\begin{thm}
\label{ThmMeanFieldRDE}
Let $V,V_1,\dots,V_\ell : \RR^d\times\mcM_1\big(\RR^d\big)\rightarrow \RR^d$ be vector fields on $\RR^d$. Given any $P\in\mcM_1\big(\RR^d\big)$, we assume that $V(\cdot,P)$ is of class $\mcC^2$, with associated norm no greated than $\lambda$, and that the vector fields $V_i(\cdot,P)$ are of class $\mcC^{[p]+1}$, with associated norms uniformly bounded \wrt $P$, and that they satisfy the inequalities
\begin{equation}
\label{EqConditionViVNonlinearRDE}
\underset{i=1..\ell}{\max}\;\big\|V_i(\cdot,P) - V_i(\cdot,Q)\big\|_\infty \vee \big\|V(\cdot,P) - V(\cdot,Q)\big\|_\infty \leq \lambda\,\textrm{\emph{d}}(P,Q), 
\end{equation}
for all $P,Q\in\mcM_1\big(\RR^d\big)$. Assume also that the random variables $\big(X^{r,I}_{ba}\big)_{r=1..[p],\,I\in\llbracket 1,\ell\rrbracket^r}$ are integrable and that there exists for each $0\leq a\leq T$, a positive random variable $C_a$, such that each of them satisfies the inequality  
\begin{equation}
\label{ConditionLawX}
\Big|\EE\Big[X^{r,I}_{ba}\,\Big|\,\mcF_a\Big]\Big|\leq C_a(b-a),
\end{equation}
for all $0\leq a\leq b\leq T$, with $\underset{0\leq a\leq T}{\sup}\,\EE[C_a]<\infty$. Last, suppose that 
\begin{equation}
\label{MomentConditionX}
\EE\big[\|\bfX\|^\gamma\big]<\infty.
\end{equation}
Then one can choose $m$ big enough so that the map $\Phi$ has a fixed point in $\textrm{\emph{Lip}}(m)$. This fixed point is unique and depends continuously on the law of $\bfX$ if the vector fields $V_i$ do not depend on their $P$-component.
\end{thm}

\ssk

\begin{Dem}
\noindent \textbf{1. Existence.} We first prove that one can choose $m$ big enough so that the map $\Phi$ sends $\textrm{Lip}(m)$ to itself. Given such a path $\mcP$, the time-dependent bounded vector field $V(x,P_t)$ is Lipschitz in its two arguments, so we can denote by $\varphi^{\mcP}$ the well-defined solution flow to the time non-homogeneous rough differential equation \eqref{NonlinearRDE} on $\RR^d$. We need to see that, for any function $g\in\mcC^\gamma\big(\RR^d\big)$, with $\|g\|_\gamma\leq 1$, and any $0\leq s\leq t\leq T$, we have
\begin{equation*}
\big(g,\Phi(\mcP)_t\big) -\big(g,\Phi(\mcP)_s\big)\leq c\,(t-s),
\end{equation*}
for some positive constant $c\leq m$, that is 
\begin{equation}
\label{LipContinuityF}
\EE\Big[g\big(\varphi_{t0}^\mcP(x_0)\big)-g\big(\varphi_{s0}^\mcP(x_0)\big)\Big] \leq c\,(t-s).
\end{equation}
Set $s_i=s+i2^{-n}(t-s)$ and define the approximate flow 
$$
\varphi_{ts}^{(n),\mcP} = \mu_{s_{2^n}s_{2^n-1}}^\mcP\circ\cdots\circ\mu_{s_1s_0}^\mcP
$$
associated with the time non-homogeneous rough differential equation \eqref{NonlinearRDE} as in remark \ref{Remarks} (6); it converges uniformly to $\varphi^\mcP$. The expectation $\EE\Big[g\Big(\varphi_{t0}^\mcP(x_0)\Big)-g\Big(\varphi_{s0}^\mcP(x_0)\Big)\Big]$ is, by dominated convergence, the limit of $\EE\Big[g\Big(\varphi_{t0}^{(n),\mcP}(x_0)\Big) - g\Big(\varphi_{s0}^{(n),\mcP}(x_0)\Big)\Big]$, and writing $g\Big(\varphi_{t0}^{(n),\mcP}(x_0)\Big)-g\Big(\varphi_{s0}^{(n),\mcP}(x_0)\Big)$ as a telescopic sum, using the elementary identity \eqref{ElementaryIdentity}, gives
\begin{equation*}
\EE\Big[g\Big(\varphi_{t0}^\mcP(x_0)\Big) - g\Big(\varphi_{s0}^\mcP(x_0)\Big)\Big] = \underset{n\rightarrow \infty}{\lim}\,\sum_{k=0}^{2^n-1}\EE\Big[g\big(\mu^\mcP_{s_{k+1}s_k}\overline\mu^\mcP_k({\bfz})\big)-g\big(\overline\mu^\mcP_k(z)\big)\Big], 
\end{equation*}
where we write $z$ for $\varphi_{s0}^\mcP(x_0)$, and $\overline\mu^\mcP_k = \mu^\mcP_{s_ks_{k-1}}\circ\cdots\circ\mu^\mcP_{s_1s_0}$, for $1\leq k\leq 2^n$. With these notations, using the Taylor expansion of $\mu^\mcP_{s_{k+1}s_k}$ given in \eqref{GralFundamentalEstimate}, we see that $\EE\big[g\big(\varphi_{t0}^\mcP(x_0)\big)-g\big(\varphi_{s0}^\mcP(x_0)\big)\big]$ is the limit of 
\begin{equation*}
\begin{split}
\sum_{k=0}^{2^n-1}\EE\Big[&2^{-n}(t-s) \big(V(\cdot,P_{s_k})g\big) \big(\overline\mu^\mcP_k(z)\big) + \\
&\sum_{r=1}^{[p]} \sum_{I\in\llbracket 1,\ell\rrbracket^r}X^{r,I}_{s_{k+1}s_k}\big(V_I(\cdot,P_{s_k})g\big)\big(\overline\mu^\mcP_k(z)\big) + \|g\|_\gamma\big(2+\|{\bfX}\|^{\gamma}\big)\,o\Big(\big(2^{-\frac{\gamma}{p}\,n}\Big)\Big],
\end{split}
\end{equation*}
where $o(\cdot)$ does not depend on $\bfX$. So, using \eqref{ConditionLawX} and \eqref{MomentConditionX} and the uniform boundedness assumption in $\mcC^{[p]+1}$ of the vector fields $V_i(\cdot,P)$, the quantity $\EE\big[g\big(\varphi_{t0}^\mcP(x_0)\big)-g\big(\varphi_{s0}^\mcP(x_0)\big)\big]$ appears as no greater than the large $n$ limit of 
\begin{equation*}
\begin{split}
&c(t-s) + \sum_{k=0}^{2^n-1}\sum_{r=1}^{[p]} \sum_{I\in\llbracket 1,\ell\rrbracket^r} \EE\Big[\EE\big[X^{r,I}_{s_{k+1}s_k}\big|\mcF_{s_k}\big]\,\big(V_I(\cdot,P_{s_k})g\big)\big(\overline\mu^\mcP_k(z)\big)\Big] + o_n(1)\\
&\leq c(t-s) + c \sum_{k=0}^{2^n-1} \EE\big[C_{s_k}\big]\,\big(s_{k+1}-s_k\big) + o_n(1) \leq c(t-s) + o_n(1).
\end{split}
\end{equation*}
 Inequality \eqref{LipContinuityF} follows by choosing $m$ big enough. Fix $P_0\in\mcM_1\big(\RR^d\big)$ and define $\textrm{Lip}_0(m)$ as the subset of $\textrm{Lip}(m)$ paths with starting point $P_0$. As $\textrm{Lip}_0(m)$ is a convex compact subset of the set of continuous paths from$[0,T]$ to $\Big(\mcM_1\big(\RR^d\big),\textrm{d}\Big)$, equipped with the norm of uniform convergence, Schauder's fixed point theorem applies and gives the existence of a fixed point of the map $\Phi$.
 
\medskip

\noindent \textbf{2. Uniqueness.} We suppose in this paragraph that the vector fields $V_i(\cdot,P)=V_i(\cdot)$ do not depend on their $\mcM_1(\RR^d)$-component. First, we prove that $\Phi$ is a strict contraction of $\textrm{Lip}_0(m)$, provided $T$ is small enough. We use for that purpose the same telescopic decomposition as above, which is reminiscent of the well-known identity $\textrm{T}_t = \textrm{Id}+ \int_0^t \textrm{A}\textrm{T}_r\,dr$, satisfied by the semi-group $({\textrm T}_t)_{t\geq 0}$ of any Markov process with generator A. 

\ssk

Fix $\omega\in\Omega$ and omit it in the notations of this paragraph. Denote by $\mu_{ts}^\mcQ$ the approximate flow associated to the time non-homogeneous rough differential equation \eqref{NonlinearRDE}, with $Q_t$ instead of $P_t$, for all $0\leq t\leq T$. Given $g\in\mcC^\gamma\big(\RR^d\big)$, we see, using dominated convergence, that
\begin{equation*} 
\begin{split}
\EE\Big[ g\big(\varphi_{t0}^\mcP(x_0)\big) - g\big(\varphi_{t0}^\mcQ(x_0)\big) \Big] = \underset{n\rightarrow \infty}{\lim}\; \EE\Big[ g\big(\varphi_{t0}^{(n),\mcP}(x_0)\big) - g\big(\varphi_{t0}^{(n),\mcQ}(x_0)\big) \Big] 
\end{split}
\end{equation*}
is the limit of
\begin{equation*} 
\begin{split}
\sum_{k=0}^{2^n-1} \EE\Big[ &\Big\{g\big(\mu_{s_{2^n}s_{2^n-1}}^\mcQ\circ\cdots\circ\mu_{s_{2^n-k+1}s_{2^n-k}}^\mcQ\circ\mu_{s_{2^n-k}s_{2n-k-1}}^\mcP\big)  \\
&- g\big(\mu_{s_{2^n}s_{2^n-1}}^\mcQ\circ\cdots\circ\mu_{s_{2^n-k+1}s_{2^n-k}}^\mcQ\circ\mu_{s_{2^n-k}s_{2^n-k-1}}^\mcQ\big)\Big\}\circ \overline\mu_{2^n-k-1}^\mcP(x_0) \Big]
\end{split}
\end{equation*}
where $\overline\mu_{2^n-k-1}^\mcP= \mu_{s_{2^n-k-1}s_{2^n-k-2}}^\mcP\circ\cdots\circ\mu_{s_1s_0}^\mcP$, with the obvious convention concerning the summand for the first and last term of the sum. It follows from the time-inhomogeneous version of proposition \ref{PropEstimatesGeneralRDE} giving the Taylor expansion of $\mu_{s_{2^n-k+1}s_{2^n-k}}^{\mcQ,\mcP}$, that $\EE\Big[ g\big(\varphi_{t0}^\mcP(x_0)\big) - g\big(\varphi_{t0}^\mcQ(x_0)\big)\Big]$ is no greater than the large $n$ limit of 
\begin{equation*} 
\EE\left[c\big(1+\|\bfX\|^{\gamma}\big)\,\|g\|_{\mcC^1}\,\sum_{k=0}^{2^n-1} \big\|V(\cdot,P_{s_{2^n-k-1}}) - V(\cdot,Q_{s_{2^n-k-1}})\big\|_\infty\,2^{-n}t \right] + c\|g\|_\gamma \,\EE\big[1+\|\bfX\|^\gamma\big]\,2^{\frac{\gamma-p}{p}\,n}.
\end{equation*}
This upper bound is no greater than
\begin{equation*} 
c\,\|g\|_{\mcC^\gamma}\,\lambda\,2^{-n}t\sum_{k=0}^{2^n-1}\textrm{d}\big(P_{s_{n-k-1}},Q_{s_{n-k-1}}\big) + o_n(1) \leq ct\,\|g\|_{\mcC^\gamma}\underset{s\in [0,t]}{\sup}\,\textrm{d}(P_s,Q_s) + o_n(1),
\end{equation*}
since $\EE\big[\|\bfX\|^{\gamma}\big]$ is finite. As a result, we have for all $t\in [0,T]$
\begin{equation*} 
\begin{split}
\textrm{d}\big(\Phi(\mcP)_t,\Phi(\mcQ)_t\big) &= \underset{\|g\|_{\mcC^\gamma}\leq 1}{\sup}\,\Big|\EE\Big[g\big(\varphi^\mcP_{t0} ({\bf x}_0)\big)\Big]-\EE\Big[g\big(\varphi^\mcQ_{t0}({\bf x}_0)\big)\Big]\Big| \\
&\leq cT\,\underset{t\in [0,T]}{\max}\textrm{d}\big(P_t,Q_t\big),
\end{split}
\end{equation*}
so $\Phi$ is a strict contraction provided $cT<1$.

\ssk

 As usual in the study of ordinary differential equations, the fact that $\Phi$ has a unique fixed point for any $T$ follows from the fact that the above condition on $T$ does not involve $P_0=Q_0$.

\ssk

As the solution to the rough differential equation \eqref{NonlinearRDE} depends continuously on $\bfX$, the map $\Phi$, considered as a function of $\mcP$ and of the parameter "law of $\bfX$", is a continuous function of its two arguments. It is an elementary result that the unique fixed point of $\Phi$ is then a continuous function of the parameter.
\end{Dem}

\ssk

\begin{rem}
So far there has been only one other work dealing with mean field stochastic rough differential equations, by Cass and Lyons \cite{CassLyons}. They prove an existence and well-posedness result under a different set of hypotheses on the law of $\bfX$ and a different topology on the set $\mcM_1\big(\RR^d\big)$. They require exponential moments for the accumulated local variation of $\bfX$, work with Wasserstein distance on $\mcM_1\big(\RR^d\big)$ and consider only mean field dynamics with a linear interaction in the drift. Their setting covers Gaussian $p$-rough paths, for $2\leq p<4$, which the above setting does not cover. On the other hand, we are able to prove an existence result for dynamics with a \emph{nonlinear} mean field interaction in both the drift and "diffusivity" vector fields. In situations where both settings apply (like Gaussian $p$-rough paths with $1\leq p\leq 2$), the assumptions of theorem \ref{ThmMeanFieldRDE} are significantly weaker than the hypotheses of \cite{CassLyons}.
\end{rem}

\end{document}